\documentclass[]{amsart}
\usepackage[utf8]{inputenc}
\usepackage{amsmath,amssymb,amsthm}
\usepackage{epsfig}
\usepackage{graphicx}
\usepackage{relsize}
\usepackage{mdframed}
\usepackage{multirow}
\usepackage[width=1.0\textwidth, font=footnotesize]{caption}
\usepackage{amsmath,amsbsy}
\usepackage{mathtools}
\usepackage{tcolorbox}
\usepackage[linesnumbered,ruled,vlined]{algorithm2e}
\SetKwInput{KwInput}{Input}                
\SetKwInput{KwOutput}{Output}              
\usepackage{hyperref}
\usepackage{subcaption}
\hypersetup{
	colorlinks=true,
	linkcolor=blue,
	citecolor=red,
	filecolor=red,      
	urlcolor=cyan,
}
\newtheorem{theorem}{Theorem}

\newtheorem*{notation}{Notation}

\theoremstyle{definition}
\newtheorem{assumption}{Assumption}
\newtheorem{remark}{Remark}
\newtheorem{definition}{Definition}

\newcommand{\norm}[1]{\left\| #1 \right\|}
\newcommand{\abs}[1]{\left| #1 \right|}

\newcommand{\prt}[1]{\left( #1 \right)}
\newcommand{\cost}[1]{\mathrm{Cost}}

\newcommand{\bN}{\mathbb{N}}
\newcommand{\bP}{\mathbb{P}}
\newcommand{\bR}{\mathbb{R}}

\newcommand{\bT}{\mathbf{T}}

\newcommand{\fell}{\ensuremath{\boldsymbol\ell}}
\newcommand{\fe}{\ensuremath{\boldsymbol e}}
\newcommand{\fone}{\ensuremath{\boldsymbol 1}}
\newcommand{\fDelta}{\ensuremath{\boldsymbol \Delta}}
\newcommand{\cB}{\mathcal{B}}
\newcommand{\cF}{\mathcal{F}}
\newcommand{\cO}{\mathcal{O}}
\newcommand{\cI}{\mathcal{I}}

\newcommand{\cN}{\mathcal{N}}

\newcommand{\vBar}{\bar v}
\newcommand{\vHat}{\hat{v}}
\newcommand{\vBarHat}{\hat{\bar{v}}}

\newcommand{\Ex}[1]{\E \left[ #1 \right]}
\newcommand{\Var}[1]{\V \left[ #1 \right]}
\newcommand{\Cov}{\overline{\mathrm{Cov}}}
\newcommand{\E}{\mathbb{E}}
\newcommand{\V}{\mathbb{V}}

\title[Multi-index Ensemble Kalman Filtering]{Multi-index Ensemble Kalman Filtering}\thanks{© 2022. This manuscript version is made available under the CC-BY-NC-ND 4.0 license https://creativecommons.org/licenses/by-nc-nd/4.0/}

\author[H. Hoel]{H{\aa}kon Hoel} \address[{H{\aa}kon Hoel}]{\newline
  Department of Mathematics, University of Oslo, Oslo, Norway
  \newline (haakonah@math.uio.no)}

\author[G. Shaimerdenova]{Gaukhar
	Shaimerdenova$^*$}\thanks{$^*$Corresponding author: G.Shaimerdenova
	(gaukhar.shaimerdenova@kaust.edu.sa)} \address[{Gaukhar
	Shaimerdenova}]{\newline Applied Mathematics and Computational
	Sciences, KAUST, Thuwal, Saudi Arabia \newline
	(gaukhar.shaimerdenova@kaust.edu.sa)}

\author[R. Tempone]{Ra\'ul Tempone} \address[{Raul Tempone}]{\newline
	Chair of Mathematics for Uncertainty Quantification, RWTH Aachen
	University, Aachen, Germany \newline (tempone@uq.rwth-aachen.de)
	\newline
	\and
	\newline
	Applied Mathematics and Computational Sciences, KAUST, Thuwal, Saudi Arabia \newline (raul.tempone@kaust.edu.sa)}

\begin{document}
	
\begin{abstract}
  In this work we combine ideas from multi-index Monte Carlo and
  ensemble Kalman filtering (EnKF) to produce a highly efficient
  filtering method called multi-index EnKF (MIEnKF).  MIEnKF is based
  on independent samples of four-coupled EnKF estimators on a
  multi-index hierarchy of resolution levels, and it may be viewed as
  an extension of the multilevel EnKF (MLEnKF) method developed by the
  same authors in 2020. Multi-index here refers to a two-index method,
  consisting of a hierarchy of EnKF estimators that are coupled in two
  degrees of freedom: time discretization and ensemble size.  Under
  certain assumptions, when strong coupling between solutions on
  neighboring numerical resolutions is attainable, the MIEnKF method
  is proven to be more tractable than EnKF and MLEnKF. Said efficiency
  gains are also verified numerically in a series of test problems.
  \bigskip
  \noindent
  
  \textbf{Key words}: Monte Carlo, multilevel, multi-index, convergence
  rates, Kalman filter, ensemble Kalman filter
  
  \noindent \textbf{AMS subject classification}: 65C30, 65Y20. 
  
\end{abstract}

\maketitle

\section{Introduction}
The ensemble Kalman filter (EnKF) is a widely used data
assimilation method for high-dimensional state-space problems with
nonlinear dynamics. Owing to its simple implementation and efficiency,
ensemble-based filtering methods have rapidly gained popularity in
geophysical sciences with applications, for example, in weather
forecasting~\cite{kalnay2003atmospheric}, atmosphere-ocean/lake
simulations~\cite{houtekamer2005atmospheric, baracchini2020data,
  hammoud2021moving}, and oil reservoir
management~\cite{aanonsen2009ensemble, ruchi2021fast}. The EnKF method
was originally proposed by
Evensen~\cite{evensen1994sequential}. Subsequently, several variants
were developed~\cite{houtekamer1998data, bishop2001adaptive,
  anderson2001ensemble}. EnKF approximates the filtering distribution
using the empirical measure of its ensemble members. The
$L^p-$convergence of the EnKF method with perturbed
observations~\cite{houtekamer1998data} has been studied in the
literature~\cite{mandel2011convergence, le2009large}.

A considerable challenge in numerical filtering methods is the
increase in simulation cost as the numerical resolution gets finer.
This challenge can be overcome by the multilevel Monte Carlo method
(MLMC)~\cite{giles2008}, which achieves substantial variance
reduction by simulating pairwise coupled realizations on a hierarchy
of temporal discretization levels. MLMC is a flexible
methodology that has been combined with many other methods and successfully
implemented in various fields: quasi-Monte
Carlo~\cite{giles2009multilevel, kuo2017, robbe2019}, sequential Monte
Carlo~\cite{beskos2017multilevel,beskos2018multilevel,moral2017multilevel,latz2018multilevel},
inverse problems and experimental design~\cite{beck2020, schaden2020,
  litvinenko2019computation, goda2020, taverniers2020}, differential
equations with randomness~\cite{jourdain2019, detommaso2019,
  barth2018, khodadadian2020, beck2020Hp, badwaik2019}, limit
theorems~\cite{alay2020, hoel2019central}, importance
sampling~\cite{hammouda2020,kebaier2018coupling, fang2020}, and
machine learning~\cite{lye2019}.

The multilevel EnKF (MLEnKF) method was introduced by Hoel et
al.~\cite{hoel2016} for stochastic differential equation models with
discrete-time observations, and an alternative version based on a
sample average of independent pairwise coupled EnKF estimators was
subsequently developed~\cite{hoel2020multilevel}. The main difference
between these two versions of MLEnKF is that~\cite{hoel2016} uses one
universal Kalman gain to update the ensemble members on all hierarchy
levels, while~\cite{hoel2020multilevel} employs one Kalman gain per
independent EnKF sample in the full MLEnKF estimator. The latter
approach introduces less correlation between all particle members of
the MLEnKF estimator, which particularly simplifies convergence
analysis and paves the way for extending MLEnKF to the multi-index
EnKF (MIEnKF) method introduced in this work. The MLEnKF method was
extended to spatiotemporal (infinite-dimensional state space)
models~\cite{chernov2020}.  Similar multilevel techniques have been 
combined with other ensemble-based filtering methods, such as
particle filters~\cite{jasra2017multilevel,ballesio2020}, transform particle
filters~\cite{gregory2016multilevel, gregory2017seamless},
multigrid~\cite{moldovan2020} and the recent extension to
the continuous-time (Kalman--Bucy) filter~\cite{chada2020multilevel}.

The successful implementation of MLMC depends on a strong pairwise
coupling between realizations on neighboring hierarchy levels, meaning a
coupling that leads to substantial variance reduction.  For
stochastic differential equations with sufficiently smooth
coefficients, this is achieved quite easily, but in more realistic
problems with low-regularity features this can be extremely
challenging if at all possible. See~\cite{fossum2020} for
multilevel data assimilation applied to reservoir
history matching, and~\cite{gao2020bi,popov2020multifidelity,
  popov2021multifidelity} for applications of MLEnKF using sampling
resolution constraints (so-called multi-fidelity
methods).

Another important method is the multi-index Monte Carlo method
(MIMC)~\cite{abdo2016}, which forms the basis of this work. MIMC
consists of a multi-index hierarchy of coupled realizations on
neighboring resolutions, and can be regarded as an extension of
MLMC. Many concepts related to particle-wise coupling in the proposed
MIEnKF method are common in the MIMC method for McKean-Vlasov
dynamics~\cite{abdo2018}.

The contributions of this work are to develop the MIEnKF method
with a subtle variance-reducing coupling idea for realizations on
neighboring resolutions, and to numerically verify the asymptotic
efficiency gains that MIEnKF achieves over EnKF and MLEnKF. The MIEnKF
method extends the recent MLEnKF method~\cite{hoel2020multilevel} by
treating not only the numerical discretization but also the EnKF
ensemble size as degrees of freedom -- resolution parameters. 
MIEnKF introduces a four-coupling of EnKF estimators (i.e., a coupling in
both degrees of freedom) that produces a stronger variance reduction
than the pairwise coupling in MLEnKF. Under certain
assumptions, MIEnKF is also shown theoretically to achieve
efficiency gains over counterparts for weak approximations of
quantities of interest (QoI) in the classic and more robust setting
of $\alpha=1$ and $\beta=2$ defined in~\cite{hoel2020multilevel},
cf.~Table~\ref{table:1}.
\begin{table}[h!]
  \centering
  \begin{tabular}{|c |c |c |c|} 
    \hline
    Methods & EnKF & MLEnKF & MIEnKF \\ [0.5ex] 
    \hline\hline
    Mean-squared error & $\cO(\epsilon^2)$& $\cO(\epsilon^2)$ &  $\cO(\epsilon^2)$\\
    \hline
    Computational cost & $\cO(\epsilon^{-3})$& $\cO(\epsilon^{-2}\abs{\log(\epsilon)}^3)$ & $\cO(\epsilon^{-2})$ \\
    \hline
  \end{tabular}
  \bigskip
  \caption{Comparison of computational costs versus errors for ensemble Kalman filtering (EnKF), multilevel EnKF (MLEnKF) and multi-index EnKF (MIEnKF) methods, cf. Section~\ref{sec:complexity}.}
  \label{table:1}
\end{table}

The rest of this work is organized as follows. In
Section~\ref{sec:problem}, the setting and notation for filtering
problem are introduced and a brief overview of the EnKF, mean-field
EnKF (MFEnKF), and MLEnKF methods is presented.
Section~\ref{sec:mienkf} describes the framework of the MIEnKF method.
Section~\ref{sec:complexity} presents theory on the performance of the
MIEnKF method, including a theorem on approximation error versus
computational cost. Section~\ref{sec:numerics} compare the performance
of MIEnKF to MLEnKF and EnKF in a series of numerical examples, and we
wrap up with concluding remarks in Section~\ref{sec:conclusion}.

\section{Problem setting}

In this section, we introduce the filtering problem of interest and
give a brief overview of relevant ensemble-based filtering methods. 

\label{sec:problem}
Let $\prt{\Omega, \cF, \bP; \{\cF_t\}_{t\geq0}}$ be a complete
probability space equipped with a filtration $\{\cF_t\}_{t\geq0}$ of
sub-$\sigma$-algebras of $\cF=\cF_{\infty}$. We denote by
$L_t^p(\Omega, \bR^k)$ the space of $\cF_t\backslash\cB^k$-measurable
functions\footnote{The function $u$ is
  $\cF_t\backslash\cB^k$-measurable iff $u^{-1}(B)\in \cF_t$ for all
  $B\in \cB^k,$ where $\cB^k$ denotes the Borel $\sigma$-algebra on
  $\bR^k$.}  $u:\Omega \rightarrow \bR^k$ with
$\Ex{\abs{u}^p}<\infty$. Given the initial value
$u_0 \in \cap_{p\geq2} L_0^p(\Omega, \bR^d)$, we consider the
discrete-time filtering problem for a system of stochastic dynamics
defined by a sequence of random maps
$\Psi_n: \bR^d \times \Omega \rightarrow \bR^d$ and observations with
additive noise:
\begin{equation*}
  \centering
  \left\{ \begin{split}
    & u_{n+1}(\omega)=\Psi_n (u_n, \omega),\\
    & y_{n+1}(\omega)= Hu_{n+1}(\omega)+\eta_{n+1},   \; \; 
  \end{split}\right.
\end{equation*}
where $\omega\in \Omega$, $n\in \bN_0:=\bN\cup \{0\}$,
$H\in \bR^{m\times d}$ is an observation operator,
$\{\eta_k\}_{k\in\bN}$ is an independent and identically distributed
(i.i.d.) sequence with $\eta_1 \sim N(0,\Gamma)$ and with the
independence property
$\{\eta_k\}_{k\in \bN} \perp \{u_k\}_{k\in \bN_0}$.
When confusion is not possible, we shall not
indicate the dependence on $\omega$ for random variables.

Let $Y_n:=(y_1, y_2, ..., y_n)$ denote the accumulated
observation data up to time $n$ using the convention that
$Y_0:=\emptyset$. The main objective of a filtering method is to track
the underlying signal $u_n$ given $Y_n$ through computing the
conditional distribution of $u_n$ give $Y_n$.
The exact filter density for this problem -- the so-called Bayes filter --
satisfies the following iterative equations:
\[
\begin{split}
  \mbox{ \textbf{Prediction}  } \quad \rho_{u_n|Y_{n-1}}(u)
  &\propto \int_{\bR^d}\rho_{u_n|u_{n-1}}(u) \rho_{u_{n-1}|Y_{n-1}}(v)dv\\
  \mbox{ \textbf{Update}  }\quad \quad \rho_{u_n|Y_n}(u)
  &\propto  \exp\big(-\big| \Gamma^{-1/2}(y_n-Hu)\big|^2/2 \big) \; \rho_{u_n|Y_{n-1}}(u).
  \end{split}
\]
We will refer to the posterior distribution in the above update step
as the true filter. In the linear-Gaussian setting, the Kalman filter
is an exact algorithm that tracks the mean and covariance of the true
filter. When $\Psi$ is nonlinear the true filter becomes non-Gaussian,
and the Kalman filter does not apply anymore. Therefore, approximation
methods are needed. Among such, particle filters converge to the true
filter in the large particle limit, but they are conjectured to
perform poorly in high dimensions~\cite{bengtsson2008curse}. The EnKF
performs more robustly than particle filters in high dimensions, but
it has poorer convergence properties.
In the large-ensemble limit, EnKF converges to the so-called mean-field
EnKF in the large-ensemble limit~\cite{le2009large,
  hoel2020multilevel}. However, due to the application of a biased 
Gaussian ansatz in the update step of
EnKF~\cite{evensen1994sequential}, the mean-field EnKF is not equal to
the true filter in nonlinear problem settings. Despite this
disparity, the EnKF is a robust and efficient method that
is popular approach among practitioners. Connections between the
mean-field EnKF and the true filter are discussed
in~\cite{law2016deterministic,hoel2020multilevel}, but there are many
open questions that remain to be studied, such as the convergence
properties of EnKF in the large-ensemble and long-time limit.

The main objective of this paper is to construct an efficient MIEnKF
method that converges weakly to the mean-field EnKF in the
large-ensemble limit. In other words, for a given QoI
$\varphi: \bR^d \rightarrow \bR$, our method approximates
\[
\E^{\bar \mu_n} [\varphi(u)]=\int_{\bR^d} \varphi(u) \bar \mu_n (du),
\]
where $\bar \mu_n$ denotes the mean-field EnKF measure at time $n$,
cf.~Section \ref{ssec:mfenkf}. 
\begin{notation}
  ~\\
  \vspace{-0.4cm}
  
  \begin{itemize}
    
  \item For $f,g:(0,\infty) \to [0,\infty)$ the notation $f \lesssim  g$ implies that there exists a
    $C>0$ such that
    \[
    f(x) \le  C g(x), \quad  \forall x \in (0,\infty).
    \]
    
  \item The notation $f \eqsim  g$ implies that $f \lesssim g$ and $g\lesssim f$.
    
  \item The expectation operator is defined by $\Ex{\cdot}$ and the variance operator (applicable to scalar-valued rv) is denoted by $\Var{\cdot}$.
  
  \item For $d \in \bN$,
    $|x|$ denotes the Euclidean norm of a vector $x \in \bR^d$.
    For $\cF\backslash  \cB^d$-measurable functions $u:\Omega \to \bR^d$ and $p\ge 1$,
    \[
    \norm{u}_{p} := \norm{u}_{L^p(\Omega, \bR^d)} = \prt{\int_{\Omega} |u(\omega)|^p \, \mathbb{P}(d\omega)}^{1/p}.
    \]
    
    
  \item  $\lceil x \rceil := \min \{ z \in \mathbb{Z} \mid z \ge x\}$.
  \end{itemize}
\end{notation}

Let $\Psi_n^N$ denote the numerical discretization of the dynamics
$\Psi_n$ using $N\geq 1$ uniform timesteps over every observation-time interval.
The following assumption ensures
that the mean-field EnKF measure $\bar \mu_n$ is well-defined
cf.~\cite[Appendix A]{hoel2020multilevel} and Section~\ref{ssec:mfenkf}:

\begin{assumption}\label{ass:Psi}
  Let $u,v \in \cap_{p\ge 2} L^p_n(\Omega, \bR^d)$ for any $n\in
  \bN_0$, $p\ge 2$, then there exists a
  constant $c_p>0$ such that for all $N \geq 1$:
  \begin{itemize} 
  \item[(i)] $
    \norm{\Psi^{N}_n(u)}_p \leq c_p (1+\norm{u}_p),
    $
    
  \item[(ii)]
    $
    \norm{\Psi^{N}_n(u)- \Psi^{N}_n(v)}_p < c_p \norm{u-v}_p .
    $
    
  \end{itemize}
\end{assumption}

\subsection{EnKF}
\label{ssec:enkf}
The EnKF method is an ensemble-based nonlinear filtering method that is an
extension of the Kalman filter. For an EnKF ensemble of size $P$, let
$v_{n,i}:=v_n(\omega_i)$ and $\vHat_{n,i}:=\vHat_n(\omega_i)$,
respectively, denote the $i-$th particle of the prediction and updated
ensemble at time $n$. Then, the EnKF algorithm with perturbed
observations and numerical dynamics $\Psi^N$ comprises the following
steps:
\begin{equation} \label{enkf:prediction}
  \centering \mbox{ \textbf{Prediction}  }
  \left\{\begin{split}
  v_{n+1,i} &= \Psi_n^N(\vHat_{n,i}), \quad i=1,2,...,P,\\
  m_{n+1} &= \frac{1}{P} \sum_{i=1}^P v_{n+1,i},\\
  C_{n+1} &=\frac{1}{P-1}\sum_{i=1}^P \prt{v_{n+1,i}-m_{n+1}} \prt{v_{n+1,i}-m_{n+1}}^\bT.
  \end{split}\right.
\end{equation}

\begin{equation}\label{enkf:update}
  \centering \mbox{ \textbf{Update}  }
  \left\{\begin{split}
  \tilde{y}_{n+1,i}&=y_{n+1}+\eta_{n+1,i}, \quad i=1,2,...,P,\\
  K_{n+1}&=C_{n+1}H^{\bT}(HC_{n+1}H^{\bT}+\Gamma)^{-1},\\
  \vHat_{n+1,i} &= (I-K_{n+1}H)v_{n+1,i}+K_{n+1}\tilde{y}_{n+1,i},\\
  \end{split}\right.
\end{equation}
where $\eta_{n+1,i}$ are i.i.d.~draws from $N(0, \Gamma)$.

The updated EnKF empirical measure is defined by 
\[
\mu_n^{N,P}(dv)=\frac{1}{P} \sum_{i=1}^{P} \delta(dv; \vHat_{n, i}),
\]
where $\delta$ is the Dirac measure centered at $\vHat_{n, i}$, and
the expectation of a QoI $\varphi: \bR^d \rightarrow \bR$ with respect
to the EnKF empirical measure is expressed as
\begin{equation}\label{enkf:measure}
  \mu_n^{N,P}[\varphi]=\frac{1}{P}\sum_{i=1}^{P} \varphi(\vHat_{n, i}).
\end{equation}
Note that $\mu_n^{N,P}[\varphi]$ is a random variable that depends on
parameters $N$ and $P$.
Under sufficient regularity, $\mu_n^{N,P}[\varphi] \rightarrow \bar \mu_n[\varphi]$
as $N,P\rightarrow \infty$,  cf.~\cite{le2009large, hoel2020multilevel, chernov2020},
where $\bar \mu_n[\varphi]$
denotes the expectation of $\varphi$ with respect to the
mean-field EnKF measure that is introduced in the next section.

\subsection{MFEnKF} 
\label{ssec:mfenkf}
The MFEnKF is the large-ensemble
and fine-discretization limit of EnKF. In the large-ensemble limit,
the Kalman gain becomes a deterministic matrix. Consequently, one may
view MFEnKF as an ensemble of i.i.d.~noninteracting particles, so
that it suffices to represent the resulting filtering distribution by
one particle. Let $\vBar_n$ and $\vBarHat_n$ denote the prediction and
updated state of a mean-field particle at time $n$, respectively.  The
following algorithm defines the MFEnKF for fully non-Gaussian models:
\begin{equation*}
  \centering \mbox{ \textbf{Prediction}  }
  \left\{\begin{split}
  \vBar_{n+1} &= \Psi_n(\vBarHat_{n})\\
  \bar{m}_{n+1} &= \Ex{\vBarHat_{n+1}},\\
  \bar{C}_{n+1} &= \Ex{\prt{\vBarHat_{n+1}-\bar{m}_{n+1}} \prt{\vBarHat_{n+1}-\bar{m}_{n+1}}^\bT}.
  \end{split}\right.
\end{equation*}

\begin{equation*}
  \centering \mbox{ \textbf{Update}  }
  \left\{\begin{split}
  \tilde{y}_{n+1}&=y_{n+1}+\tilde{\eta}_{n+1}, \\
  \bar{K}_{n+1}&=\bar{C}_{n+1}H^{\bT}(H\bar{C}_{n+1}H^{\bT}+\Gamma)^{-1},\\
  \vBarHat_{n+1,i} &= (I-\bar{K}_{n+1}H)\vBarHat_{n+1,i}+\bar{K}_{n+1}\tilde{y}_{n+1,i},\\
  \end{split}\right.
\end{equation*}
where $\tilde{\eta}_{n+1}$ is i.i.d.~draws from $N(0, \Gamma)$.

The expectation of a QoI $\varphi:\bR^d \rightarrow \bR$ with respect to the updated mean-field EnKF measure is given by
\[
\bar\mu_n[\varphi]:=\E^{\bar{\mu}_n}[\varphi(v)]=\int_{\bR^d}\varphi(v) \bar\mu_n(dv).
\]

\begin{remark}
  For the EnKF filter with the numerical-solution dynamics $\Psi^N_n$
  (instead of $\Psi$), Assumption~\ref{ass:Psi} ensures that the
  analogous mean-field EnKF measure $\bar \mu_n^{N}$ is well-defined
  for any observation time $n\ge 0$ and numerical resolution $N\ge 1$,
  cf.~\cite[Appendix A]{hoel2020multilevel}.
\end{remark}

\subsection{MLEnKF} \label{ssec:mlenkf} The recently developed MLEnKF
method~\cite{hoel2020multilevel} is a natural stepping stone on the
way from EnKF to explaining all of the complexities in the MIEnKF method.
MLEnKF is a filtering method based on a sample average of independent and
pairwise coupled samples of EnKF estimators at different resolution
levels.

Let $L\in \bN$ denote the finest
resolution level of the estimator, and let the sequences
\[
N_\ell = N_0\times 2^\ell 
\quad \text{and} \quad P_\ell = P_0 \times 2^\ell \quad \text{with} \quad N_0,P_0 \in \bN, \quad
\ell=0,1, \ldots ,L
\]
respectively denote the numerical resolution and ensemble size.

\subsubsection*{Pairwise coupling of EnKF estimators} For a level $\ell\ge 0$,
let 
\[
\vHat_{n,i}^\ell:=\vHat_n^\ell(\omega_i^\ell)  \qquad i =1, \ldots, P_\ell
\]
denote $i-$th particle of the updated ensemble at time $n$ in size
$P_\ell$ corresponding to the fine-level numerical resolution
$N_\ell$. Each $\vHat_{n,i}^\ell$ is \textit{coupled pairwise} to the
respective $i-$th particle of the coarser-level updated ensemble at
time $n$ computed with the numerical resolution $N_{\ell-1}$ via
shared driving noise $\omega_i^\ell$. To obtain a
$1 \leftrightarrow 1$ coupling between ensemble-members/particles on the fine-
and coarse level, the total size of the coarse-level ensemble
is set to $P_\ell$ with the relation $P_\ell=2P_{\ell-1}$, meaning
that the coarse-level ensemble can be viewed as a union of two
ensembles in size $P_{\ell-1}$:
\[
\vHat_{n,i}^{\ell-1,1}:= \vHat_n^{\ell-1,1}(\omega_i^\ell)  \qquad i =1, \ldots, P_{\ell-1} 
\]
and
\[
\vHat_{n,i}^{\ell-1,2}:= \vHat_n^{\ell-1,2}(\omega_{P_{\ell-1}+i}^\ell)  \qquad i =1, \ldots, P_{\ell-1},
\]
with the convention $\vHat^{-1,\cdot}:=0$.

It is important to note here that the particle-wise pairs share the
same realization of driving noise within a level and the superscript
$\ell$ in $\omega_i^\ell$ indicates an independence of underlying
noise between levels. In addition to this, the pairwise coupling is
imposed under the particle-wisely shared initial condition:
\[
\vHat_{0,i}^\ell = \begin{cases} \vHat_{0,i}^{\ell-1,1} & \text{if} \quad  i \in\{ 1, \ldots, P_{\ell-1}\}\\
  \vHat_{0,i-P_{\ell-1}}^{\ell-1,2} & \text{if} \quad i \in\{ P_{\ell-1}+ 1, \ldots, P_{\ell}\},
  \end{cases}
\]
and the perturbed observations are also shared particle-wisely (see the below update step).
Iterative simulation of pairwise coupled ensemble-members on the $\ell$-th resolution level of the MLEnKF filter
consists of the following prediction and update steps: 
\begin{equation}\label{mlenkf:prediction}
  \centering \mbox{ \textbf{Prediction}  }
  \left\{\begin{split}
  v_{n+1,i}^{\ell-1,1} &= \Psi_n^{N_{\ell-1}}(\vHat_{n,i}^{\ell-1,1}),  \qquad i=1,\ldots,P_{\ell-1},\\
  v_{n+1,i}^{\ell-1,2} &= \Psi_n^{N_{\ell-1}}(\vHat_{n,i}^{\ell-1,2}),  \qquad i=1,\ldots,P_{\ell-1},\\
  v_{n+1,i}^{\ell} &= \Psi_n^{N_\ell}(\vHat_{n,i}^{\ell}),\qquad \qquad   i=1,\ldots ,P_\ell, \\
  C_{n+1}^{\ell-1,1} &=\Cov[v_{n+1, 1:P_{\ell-1}}^{\ell-1,1}],\\
  C_{n+1}^{\ell-1,2} &=\Cov[v_{n+1,1:P_{\ell-1}}^{\ell-1,2}],\\
  C_{n+1}^{\ell} &=\Cov[v_{n+1,1:P_\ell}^{\ell}],
  \end{split}\right.
\end{equation}

\begin{equation*}
  \begin{split}
    \Cov[v_{n,1:P_\ell}^{\ell}]&:=\sum_{i=1}^{P_\ell}\frac{(v_{n,i}^{\ell})(v_{n,i}^{\ell})^\bT}{P_\ell}-\prt{\sum_{i=1}^{P_\ell}\frac{v_{n,i}^\ell}{P_\ell}}\prt{\sum_{i=1}^{P_\ell}\frac{v_{n,i}^{\ell}}{P_\ell}}^\bT,\\
    \Cov[v_{n,1:P_{\ell-1}}^{\ell-1,k}]&:=\sum_{i=1}^{P_{\ell-1}}\frac{(v_{n,i}^{\ell-1,k})(v_{n,i}^{\ell-1,k})^\bT}{P_{\ell-1}}-\prt{\sum_{i=1}^{P_{\ell-1}}\frac{v_{n,i}^{\ell-1,k}}{P_{\ell-1}}}\prt{\sum_{i=1}^{P_{\ell-1}}\frac{v_{n,i}^{\ell-1,k}}{P_{\ell-1}}}^\bT, \quad k=1,2.
  \end{split}
\end{equation*}
\begin{equation}\label{mlenkf:update}
  \centering \mbox{ \textbf{Update} }
  \left\{\begin{split}
  \tilde{y}_{n+1,i}^\ell&=y_{n+1}+\eta_{n+1,i}^\ell, \qquad \qquad \qquad \qquad \qquad \qquad \quad  i=1,\ldots,P_\ell, \\
  K_{n+1}^{\ell-1,1}&=C_{n+1}^{\ell-1,1}H^{\bT}(HC_{n+1}^{\ell-1,1} H^{\bT}+\Gamma)^{-1},\\
  K_{n+1}^{\ell-1,2}&=C_{n+1}^{\ell-1,2} H^{\bT}(HC_{n+1}^{\ell-1,2} H^{\bT}+\Gamma)^{-1},\\
  K_{n+1}^\ell&=C_{n+1}^\ell H^{\bT}(HC_{n+1}^\ell H^{\bT}+\Gamma)^{-1},\\
  \vHat_{n+1,i}^{\ell-1,1} &= (I-K_{n+1}^{\ell-1,1}H)v_{n+1,i}^{\ell-1,1}+K_{n+1}^{\ell-1,1}\tilde{y}_{n+1,i}^\ell,\quad  \qquad i=1,\ldots,P_{\ell-1},\\
  \vHat_{n+1,i}^{\ell-1,2} &= (I-K_{n+1}^{\ell-1,2}H)v_{n+1,i}^{\ell-1,2}+K_{n+1}^{\ell-1,2}\tilde{y}_{n+1,P_{\ell-1}+i}^\ell, \;\; i=1,\ldots,P_{\ell-1},\\
  \vHat_{n+1,i}^{\ell} &= (I-K_{n+1}^{\ell}H)v_{n+1,i}^{\ell}+K_{n+1}^{\ell}\tilde{y}_{n+1,i}^\ell,\qquad \quad \quad  i=1,\ldots,P_\ell,
  \end{split}\right.
\end{equation}
where $\{\eta_{n+1,i}^\ell\}_{i=1}^{P_\ell}$ is a sequence of independent $N(0, \Gamma)$-distributed random variables.
In the above notation, the coupling between the fine-level EnKF estimator
\[
\mu_n^{N_\ell, P_\ell}[\varphi] := \sum_{i=1}^{P_\ell} \frac{\varphi(\vHat_{n,i}^{\ell})}{P_\ell},
\]
and the two coarse-level estimators
\[
\mu_n^{N_{\ell-1}, P_{\ell-1},k}[\varphi] := \sum_{i=1}^{P_{\ell-1}} \frac{\varphi(\vHat_{n,i}^{\ell-1,k})}{P_{\ell-1}}, \qquad k=1,2
\]
is obtained through particle-wise coupling
\[
\vHat_{n,i}^{\ell} \xleftrightarrow{\text{coupling}} \begin{cases} \vHat_{n,i}^{\ell-1,1} & \text{if} \quad  i \in\{ 1, \ldots, P_{\ell-1}\}, \\
    \vHat_{n,i-P_{\ell-1}}^{\ell-1,2} & \text{if} \quad i \in\{ P_{\ell-1}+ 1, \ldots, P_{\ell}\}.
  \end{cases}
\]
Finally, the updated MLEnKF estimator at time $n$ assumes the following form:
\begin{equation}\label{eq:updMLEnKF}
	\mu_n^{ML}[\varphi] = \sum_{\ell=0}^L \sum_{m=1}^{M_\ell} \frac{\Big( \mu_n^{N_\ell, P_\ell,m} - (\mu_n^{N_{\ell-1}, P_{\ell-1},1,m} + \mu_n^{N_{\ell-1}, P_{\ell-1},2,m})/2 \Big)[\varphi]}{M_\ell}
\end{equation}
where $\{M_\ell\}_{\ell=0}^L \subset \bN$ is a decreasing sequence with $M_\ell$ representing the number of i.i.d.~and pairwise coupled EnKF estimators
on level $\ell$: 
\[
\Big\{ \mu_n^{N_\ell, P_\ell,m}[\varphi] , (\mu_n^{N_{\ell-1}, P_{\ell-1},1,m}[\varphi], \mu_n^{N_{\ell-1}, P_{\ell-1},2,m}[\varphi]) \Big\}_{m=1}^{M_\ell},
\]
where $\mu_n^{N_{-1}, P_{-1},m}[\varphi] := 0$. For the configuration of $L$ and $M_\ell$, we refer the reader to~\cite[Corollary 2]{hoel2020multilevel}. 

\section{MIEnKF}
\label{sec:mienkf}
In this section, we develop the MIEnKF method by extending
the MLEnKF method from the previous section. 

To define a set of discretization levels for the MIEnKF, we
first introduce the 2-index $\fell:=(\ell_1, \ell_2)\in \bN_0^2$
with the shorthands $\fe_1:=(1,0)$, $\fe_2:=(0,1)$, and
$\fone:=(1,1)$. Similarly as for MLEnKF, we associate sequences of
natural numbers $N_{\ell_1} = N_0\times2^{\ell_1}$ and $P_{\ell_2} = P_0\times
2^{\ell_2}$ with $N_0,P_0 \in \bN$ to the number of timesteps and particles on the 2-index ``level''
$\fell$.

Seeking to approximate $\bar{\mu}_n[\varphi]$ (the expectation of the
QoI $\varphi$ with respect to the mean-field measure $\bar \mu_n$),
we denote the discrete approximation corresponding to the 2-index
$\fell$ by
$\mu_n^{\fell}[\varphi]:=\mu_n^{N_{\ell_1}, P_{\ell_2}}[\varphi]$. In
other words, $\mu_n^{\fell}[\varphi]$ is the EnKF
estimator~\eqref{enkf:measure} computed with $N_{\ell_1}$ timesteps
and $P_{\ell_2}$ ensemble-members/particles. We define first-order difference operators
for numbers of timesteps and particles as follows:
\begin{equation}\label{eq:differenceOperators}
  \begin{split}
    \Delta_1\mu_n^{\fell}[\varphi]&=
    \begin{cases}
      \prt{\mu_n^{\fell}-\mu_n^{\fell-\fe_1}}[\varphi], \qquad \qquad \qquad \quad &\mbox{  if } \ell_1>0, \\
      \mu_n^{\fell}[\varphi], &\mbox{  if } \ell_1=0
    \end{cases}\\
    \Delta_2\mu_n^{\fell}[\varphi]&=
    \begin{cases}
      \prt{\mu_n^{\fell}-\prt{\mu_n^{\fell-\fe_2,1}+\mu_n^{\fell-\fe_2,2}}/2}[\varphi], &\mbox{  if } \ell_2>0, \\
      \mu_n^{\fell}[\varphi], &\mbox{  if } \ell_2=0
    \end{cases}
  \end{split}
\end{equation}
where $\mu_n^{\fell-\fe_2,1}[\varphi]$ and
$\mu_n^{\fell-\fe_2,2}[\varphi]$ are two i.i.d.~copies of
$\mu_n^{\fell-\fe_2}[\varphi]$. Note that $\mu_n^{\fell}[\varphi]$
comprises $P_{\ell_2}$ ensemble members, whereas
$\mu_n^{\fell-\fe_2}[\varphi]$ has half as many ensemble members, $P_{\ell_2-1}$.
Therefore, the pair $(\mu_n^{\fell-\fe_2,1}[\varphi],\mu_n^{\fell-\fe_2,2}[\varphi])$
are introduced to achieve a $1\leftrightarrow 1$ coupling of ensemble-members/particles on
the ``levels''/2-indices $\fell$ and $\fell - \fe_2$. 

We define the four-coupled EnKF estimator using the first-order mixed difference as follows:
\begin{equation}\label{eq:mixDiffOper}
  \begin{split}
    \fDelta \mu_n^{\fell}[\varphi]&:=\Delta_1(\Delta_2 \mu_n^{\fell}[\varphi])=\Delta_2(\Delta_1 \mu_n^{\fell}[\varphi])=\Delta_2(\mu_n^{\fell}-\mu_n^{\fell-\fe_1})[\varphi]\\
    &=\Bigg(\mu_n^{\fell} - \Big(\mu_n^{\fell-\fe_2,1} + \mu_n^{\fell-\fe_2,2}\Big)/2\\
    &\qquad - \mu_n^{\fell-\fe_1} + \Big(\mu_n^{\fell-\fone,1} + \mu_n^{\fell-\fone,2}\Big)/2 \Bigg)[\varphi],
  \end{split}
\end{equation}
where the pair $(\mu_n^{\fell-\fone,1}[\varphi],
\mu_n^{\fell-\fone,2}[\varphi])$ of i.i.d.~copies of
$\mu_n^{\fell-\fone}[\varphi]$  is also introduced to achieve
 a $1\leftrightarrow 1$ coupling of ensemble-members/particles on
the ``levels''/2-indices $\fell$ and $\fell - \fone$.
If it holds that 
\[
  \Ex{\mu_n^{\fell}[\varphi]} \rightarrow \bar{\mu}_n[\varphi] \quad \text{ as } \quad (\ell_1, \ell_2)  \to (\infty, \infty),
\]
and the sequence $\{ \Ex{\fDelta \mu_n^{\fell}[\varphi]} \}_{\fell \in \bN_0^2}$ is absolutely summable 
(both of these conditions hold under Assumption~\ref{ass:Psi2} ~\eqref{ass:A1}, which is presented below),
then the linearity of the expectation operator implies that 
\begin{equation}\label{eq:telescsum}
  \bar{\mu}_n[\varphi]=\sum_{ \fell  \in \bN_0^2} \Ex{\fDelta \mu_n^{\fell}[\varphi]} =\sum_{ \fell  \in \cI} \Ex{\fDelta \mu_n^{\fell}[\varphi]}+\sum_{ \fell  \notin \cI} \Ex{\fDelta \mu_n^{\fell}[\varphi]},
\end{equation}
for any index set $\cI\subset \bN_0^2$.

\begin{remark}
  The magnitude of the second term on the right-hand side, the
  truncated region, relates to the bias error of the MIEnKF estimator.
  Accurate information on how this magnitude varies with $\cI$ can and
  should be used to determine an index set such that the resulting
  MIEnKF method meets the bias-error constraint. In the last part of
  the proof of Theorem~\ref{thm:complexityMIEnKF} below, we have
  indeed used Assumption~\ref{ass:Psi2} ~\eqref{ass:A1} to control the
  bias error through bounding the magnitude of said right-hand-side term 
  in~\eqref{eq:telescsum}.
\end{remark}

For a given index set
$\cI$, which we will specify later in Section~\ref{sec:complexity}, the
MIEnKF estimator is defined as the sample-average estimator of the
first term on the right-hand side of~\eqref{eq:telescsum}:
\begin{equation} \label{eq:MIestimator}
    \mu^{MI}_n[\varphi] := \sum_{ \fell  \in \cI}
    \sum_{m=1}^{M_{\fell}} \frac{\fDelta \mu_n^{\fell,m}[\varphi]}{M_{\fell}},
\end{equation}
where $\{\fDelta \mu_n^{\fell,m}[\varphi]\}_{m=1}^{M_{\fell}}$ are i.i.d. copies of $\fDelta \mu_n^{\fell,m}[\varphi]$,
and 
$\{ \fDelta \mu_n^{\fell,m}[\varphi]\}_{(\fell,m)}$ are mutually independent.

The primary motivation for sampling four-coupled EnKF estimators in
the MIEnKF estimator is that it leads to a substantial variance
reduction that improves the tractability of the sampling method.
Similar to multilevel Monte Carlo estimators, the tractability
of~\eqref{eq:MIestimator} is optimized through careful selection of
the index set $\cI$ and the number of samples $M_{\fell}$. Provided that 
convergence rates for the multi-index hierarchy are available or
approximable, this can be achieved by solving a constrained optimization
problem~\cite{abdo2016,abdo2018}.

\subsection{Four-coupled EnKF estimators}\label{subsec:fourCoupledEstimators}
To describe the coupling between the EnKF estimators 
$(\mu_n^{\fell}, \mu_n^{\fell-\fe_2},
\mu_n^{\fell-\fe_1},\mu_n^{\fell-\fone})[\varphi]$, we
introduce the four-coupled updated-state ensembles
at time $n$:
\[
  \{(\vHat_{n}^{\fell},
\vHat_{n}^{\fell-\fe_2},\vHat_{n}^{\fell-\fe_1},\vHat_{n}^{\fell-\fone})_i\}_{i=1}^{P_{\ell_2}}:=\{(\vHat_{n}^{\fell},
\vHat_{n}^{\fell-\fe_2},\vHat_{n}^{\fell-\fe_1},\vHat_{n}^{\fell-\fone})(\omega_i^{\fell})\}_{i=1}^{P_{\ell_2}}.
\]

The set of particles on index $\fell-\fe_2$ is a union of two EnKF ensembles:
\[
\vHat_{n, i}^{\fell-\fe_2,1} := \vHat_{n, i}^{\fell-\fe_2}  \qquad i = 1, \ldots P_{\ell_2-1},
\]
and
\[
\vHat_{n, i}^{\fell-\fe_2,2} := \vHat_{n, P_{\ell_2-1} + i}^{\fell-\fe_2}  \qquad i = 1, \ldots P_{\ell_2-1},
\]
and the set of particles on index $\fell - \fone$ is also a union of two EnKF ensembles:
\[
\vHat_{n, i}^{\fell-\fone,1} := \vHat_{n, i}^{\fell-\fone}  \qquad i = 1, \ldots P_{\ell_2-1},
\]
and
\[
\vHat_{n, i}^{\fell-\fone,2} := \vHat_{n, P_{\ell_2-1} + i}^{\fell-\fone}  \qquad i = 1, \ldots P_{\ell_2-1}.
\]

Similarly to MLEnKF (Section~\ref{ssec:mlenkf}), 
the MIEnKF employs a $1 \leftrightarrow 1$ coupling between particles on all four levels
associated to one index $\fell$.  We defer further details on how this
is achieved to Section~\ref{ssec:mienkfalg}, and are now ready to
properly define the MIEnKF estimator.


The empirical estimator $\mu_n^{\fell} [\varphi]$ is induced by the
ensemble $\vHat_{n, 1:P_{\ell_2}}^{\fell} :=
\{\vHat_{n,i}^{\fell}\}_{i=1}^{P_{\ell_2}}$, meaning that it equals the
sample average of
$\{\varphi(\vHat_{n,i}^{\fell})\}_{i=1}^{P_{\ell_2}}$.  Similarly,
$\mu_n^{\fell-\fe_1}[\varphi]$ is induced by $\vHat_{n,
  1:P_{\ell_2}}^{\fell-\fe_1}:= \{\vHat_{n,
  i}^{\fell-\fe_1}\}_{i=1}^{P_{\ell_2}}$, 
\[
\mu_n^{\fell-\fe_2}[\varphi]:=\frac{(\mu_n^{\fell-\fe_2,1} +
  \mu_n^{\fell-\fe_2,2})[\varphi]}{2}
\]
is induced by the union of two
ensembles
\[\vHat_{n, 1:P_{\ell_2}}^{\fell-\fe_2}:=\{\vHat_{n,
  k}^{\fell-\fe_2,1}\}_{k=1}^{P_{\ell_2-1}}\cup \{\vHat_{n, k}^{\fell-\fe_2,2}\}_{k=1}^{P_{\ell_2-1}},
\]
and
\[
\mu_n^{\fell-\fone}[\varphi]:=\frac{(\mu_n^{\fell-\fone,1} + \mu_n^{\fell-\fone,2})[\varphi]}{2}
\]
is induced by $\vHat_{n, 1:P_{\ell_2}}^{\fell-\fone}:=\{\vHat_{n,
  k}^{\fell-\fone,1}\}_{k=1}^{P_{\ell_2-1}}\cup \{\vHat_{n,
  k}^{\fell-\fone,2}\}_{k=1}^{P_{\ell_2-1}}$.
For consistency with~\eqref{eq:differenceOperators}, we impose the condition that 
$\mu_n^{\fell-\fe_1}[\varphi]=\mu_n^{\fell-\fone}[\varphi]=0$ when
$\ell_1=0$, and $\mu_n^{\fell-\fe_2}[\varphi]=0$ when $\ell_2=0$. Figure~\ref{fig:mienkfDiag} shows a visual description
of all the couplings of the MIEnKF estimator. 

Then, the MIEnKF estimator~\eqref{eq:MIestimator} can also be written
\begin{equation}\label{eq:MIestimshort}
  \begin{split}
    \mu^{MI}_n[\varphi] := \sum_{ \fell \in \cI}
    \sum_{m=1}^{M_{\fell}} \frac{\big(\mu_n^{\fell,m}-\mu_n^{\fell-\fe_1,m}-\mu_n^{\fell-\fe_2,m}+\mu_n^{\fell-\fone,m}\big)[\varphi]}{M_{\fell}}
  \end{split}
\end{equation}
where $\{(\mu_n^{\fell,m},\mu_n^{\fell-\fe_1,m},\mu_n^{\fell-\fe_2,m},\mu_n^{\fell-\fone,m})[\varphi]\}_m$
are independent copies of the estimators $(\mu_n^{\fell},\mu_n^{\fell-\fe_1},\mu_n^{\fell-\fe_2},\mu_n^{\fell-\fone})[\varphi]$
and
$\{(\mu_n^{\fell,m},\mu_n^{\fell-\fe_1,m},\mu_n^{\fell-\fe_2,m},\mu_n^{\fell-\fone,m})[\varphi]\}_{(\fell, m)}$ are mutually independent.

\begin{remark}
  For comparison, the MLEnKF estimator~\eqref{eq:updMLEnKF} takes the following form when represented in the above 2-index notation
  \[
  \mu_n^{ML}[\varphi]=\sum_{\ell=0}^{L} \sum_{m=1}^{M_{\ell}} \frac{\big(\mu_n^{(\ell, \ell),m}-\mu_n^{(\ell,\ell) - \fone,m}\big)[\varphi]}{M_{\ell}}.
  \]
  
\end{remark}

\subsection{Particle-wise four-coupling for MIEnKF}\label{ssec:mienkfalg}
We now describe how the four-coupling of EnKF estimators manifests itself
particle by particle. At time $n=0$, the fine-index update
ensemble $\{\hat{v}^{\fell }_{0,i}\}_{i=1}^{P_{\ell_2}}$ comprises
independent $\mathbb{P}_{u_0|Y_0}$-distributed samples that are
particle-wisely coupled to three other ensembles by
\[
\hat{v}^{\fell}_{0,i}=\hat{v}^{\fell-\fe_1 }_{0,i}=\hat{v}^{\fell-\fe_2
}_{0,i}=\hat{v}^{\fell-\fone }_{0,i} \qquad \text{for} \quad  i
=1,2\ldots,P_{\ell_2}.
\]
To describe how the coupling enters in the prediction-update iterations
of the particles, let us consider the update state of a foursome
$(\vHat_{n, i}^{\fell}, \hat{v}^{\fell-\fe_1 }_{n,i}, \vHat_{n,i}^{\fell-\fe_2},\hat{v}^{\fell-\fone}_{n,i})$ at time $n\ge0$.
The next-time prediction state of the foursome is given by 
\begin{equation}
  \begin{split}
    v^{\fell }_{n+1,i}&= \Psi^{N_{\ell_1}}_n(\hat{v}^{\fell }_{n,i}), \qquad
    v^{\fell-\fe_1}_{n+1,i}= \Psi^{N_{\ell_1-1}}_n(\hat{v}^{\fell-\fe_1 }_{n,i}), \\
    v^{\fell-\fe_2}_{n+1,i} &= \Psi^{N_{\ell_1}}_n (\hat{v}^{\fell-\fe_2}_{n,i}), \qquad
    v^{\fell-\fone}_{n+1,i} = \Psi^{N_{\ell_1-1}}_n (\hat{v}^{\fell-\fone}_{n,i}),
  \end{split}
\label{mi:prediction1}
\end{equation}
where the four particles share the same driving noise in the dynamics.
The sample covariance matrices and Kalman gains are expressed as follows:
\begin{equation}
  \begin{split}
    C_{n+1}^{\fell}&=\Cov [v^{\fell}_{n+1}], \qquad K_{n+1}^{\fell }=C_{n+1}^{\fell}H^\bT(HC_{n+1}^{\fell}H^\bT+\Gamma)^{-1},\\
    C_{n+1}^{\fell-\fe_1 }&=\Cov[v^{\fell-\fe_1}_{n+1}], \! \qquad K_{n+1}^{\fell-\fe_1 }= C_{n+1}^{\fell-\fe_1  }H^\bT(HC_{n+1}^{\fell-\fe_1 }H^\bT+\Gamma)^{-1},
  \end{split}
\label{mi:prediction2}
\end{equation}
and
\begin{equation}
  \begin{split}
    C_{n+1}^{\fell-\fe_2,1}&=\Cov [v^{\fell-\fe_2,1}_{n+1}], \qquad K_{n+1}^{\fell-\fe_2,1}=C_{n+1}^{\fell-\fe_2,1}H^\bT(HC_{n+1}^{\fell-\fe_2,1}H^\bT+\Gamma)^{-1},\\
    C_{n+1}^{\fell-\fe_2,2}&=\Cov [v^{\fell-\fe_2,2}_{n+1}], \qquad K_{n+1}^{\fell-\fe_2,2}=C_{n+1}^{\fell-\fe_2,2}H^\bT(HC_{n+1}^{\fell-\fe_2,2}H^\bT+\Gamma)^{-1},\\
    C_{n+1}^{\fell-\fone,1}&=\Cov[v^{\fell-\fone,1}_{n+1}], \qquad \quad K_{n+1}^{\fell-\fone,1}= C_{n+1}^{\fell-\fone,1}H^\bT(HC_{n+1}^{\fell-\fone,1}H^\bT+\Gamma)^{-1},\\
    C_{n+1}^{\fell-\fone,2}&=\Cov[v^{\fell-\fone,2}_{n+1}], \qquad  \quad K_{n+1}^{\fell-\fone,2}= C_{n+1}^{\fell-\fone,2}H^\bT(HC_{n+1}^{\fell-\fone,2}H^\bT+\Gamma)^{-1},
  \end{split}
\label{mi:prediction3}
\end{equation}
where we recall that index $\fell-\fe_2$ and index $\fell -\fone$ both consist of two EnKF ensembles of size $P_{\ell_2-1}$,\
cf.~Section~\ref{subsec:fourCoupledEstimators}.
The perturbed observations are also particle-wisely coupled, so that one obtains the updated states: 
\begin{equation}
  \left.\begin{split}
    \tilde{y}_{n+1,i}^{\fell}&=y_{n+1}+\eta_{n+1,i}^{\fell}, \\
    \vHat_{n+1,i}^{\fell}&=(I-K_{n+1}^{\fell }H)v^{\fell}_{n+1,i}+K_{n+1}^{\fell }\tilde{y}_{n+1,i}^{\fell},\\
    \vHat_{n+1,i}^{\fell-\fe_1}&=(I-K_{n+1}^{\fell-\fe_1}H)v_{n+1,i}^{\fell-\fe_1}+K_{n+1}^{\fell-\fe_1}\tilde{y}_{n+1,i}^{\fell},
  \end{split}\right\} \quad i =1,\ldots, P_{\ell_2} \, ,
  \label{mi:update1}
\end{equation}
and
\begin{equation}
  \left.\begin{split}
    \vHat_{n+1,i}^{\fell-\fe_2,1}&=(I-K_{n+1}^{\fell-\fe_2,1 }H)v^{\fell-\fe_2,1}_{n+1,i}+K_{n+1}^{\fell-\fe_2,1 }\tilde{y}_{n+1,i}^{\fell},\\
    \vHat_{n+1,i}^{\fell-\fe_2,2}&=(I-K_{n+1}^{\fell-\fe_2,2}H)v^{\fell-\fe_2,2}_{n+1,i}+K_{n+1}^{\fell-\fe_2,2 }\tilde{y}_{n+1,i+P_{\ell_2-1}}^{\fell},\\
    \vHat_{n+1,i}^{\fell-\fone,1}&=(I-K_{n+1}^{\fell-\fone,1}H)v_{n+1,i}^{\fell-\fone,1}+K_{n+1}^{\fell-\fone,1}\tilde{y}_{n+1,i}^{\fell},\\
    \vHat_{n+1,i}^{\fell-\fone,2}&=(I-K_{n+1}^{\fell-\fone,2}H)v_{n+1,i}^{\fell-\fone,2}+K_{n+1}^{\fell-\fone,2}\tilde{y}_{n+1,i+P_{\ell_2-1}}^{\fell},
  \end{split}\right\} \quad i =1,\ldots, P_{\ell_2-1} \, ,
  \label{mi:update2}
\end{equation}
where $\{\eta_{n+1,i}^{\ell_2}\}_{i=1}^{P_{\ell_2}}$ are i.i.d.~with
$\eta_{n+1,1}^{\ell_2} \sim N(0,\Gamma)$.

\begin{figure}[htp]
  \centering
  \includegraphics[width=1.0\textwidth]{{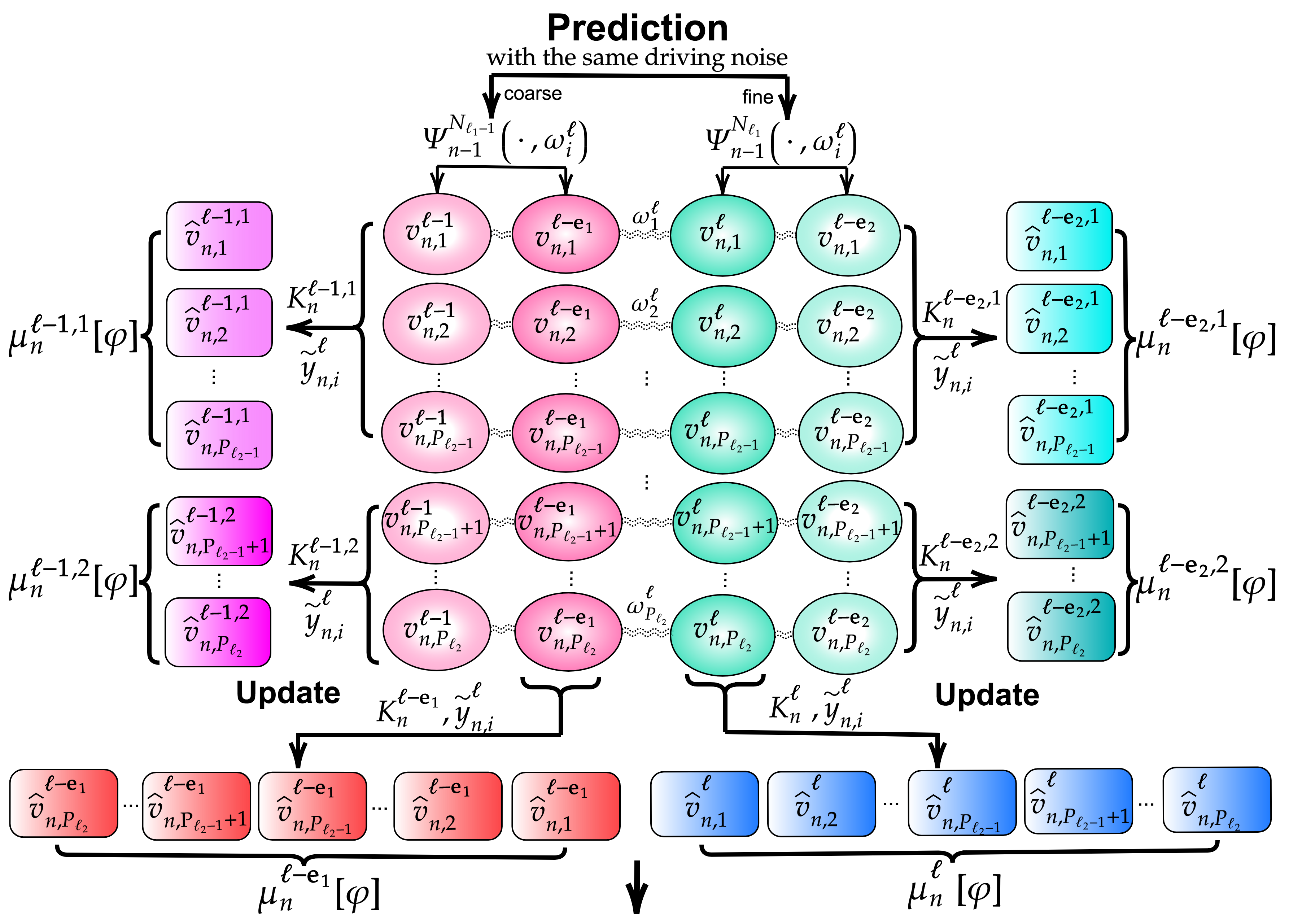}}\\
  \includegraphics[width=1.0\textwidth]{{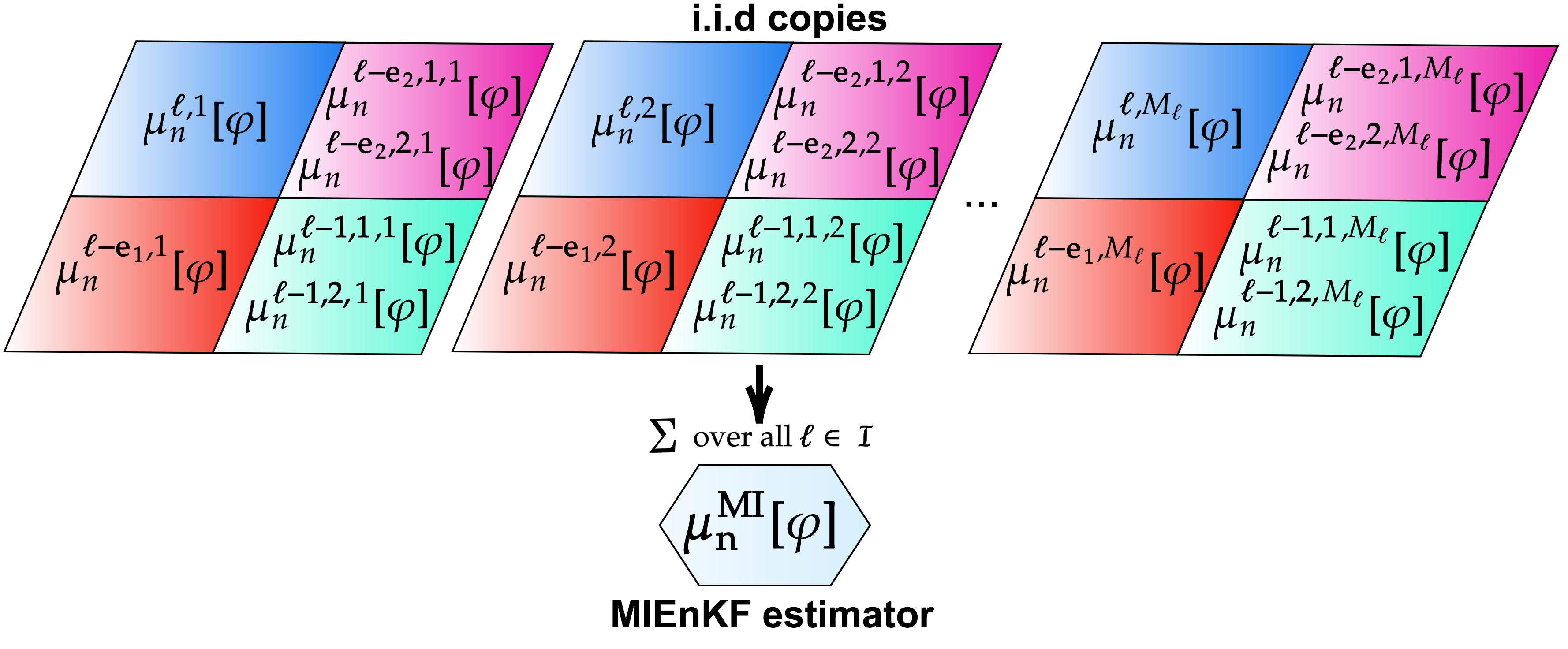}}
  \caption{One prediction-update iteration of the multi-index ensemble Kalman filtering (MIEnKF) estimator described in Section~\ref{ssec:mienkfalg}. The ovals represent four-coupled prediction-state particles, sharing the same driving noise $\omega^{\fell}$ and coupled initial conditions. The respective squares represent updated-state particles sharing the same perturbed observations.
  }
  \label{fig:mienkfDiag}
\end{figure}

To summarize, four ensembles are particle-wisely coupled by sharing the initial
condition, the driving noise, and perturbed observations. A sketch of
one prediction-update iteration of the MIEnKF method and the
composition of the MIEnKF estimator is provided in
Figure~\ref{fig:mienkfDiag} and Algorithm~\ref{alg:MIEnKF}
describes the essential steps of the MIEnKF method.

\begin{algorithm}[htp]	
	\caption{MIEnKF}\label{alg:MIEnKF}
	\DontPrintSemicolon
	
	\KwInput{The model parameters, the QoI $\varphi$, the final time $\cN$, the observation operator $H$, observations $\{y_n\}_{n=1}^{\cN}$, $m_0,\Sigma_0$, $\Gamma$, $L$, $M_{\fell}$, $N_{\ell_1}$, $P_{\ell_2}$.}
	\KwOutput{The MIEnKF estimator $\mu_n^{MI}[\varphi]$.}
	\For{$\fell \in \cI$}{
		\For{$m=1:M_{\fell}$}{
			Initialize the ensembles at time $n=0$ by sampling $\hat{v}^{\fell,m}_{0,i} \sim N(m_0,\Sigma_0)$ and setting
			$\hat{v}^{\fell,m}_{0,i}=\hat{v}^{\fell-\fe_1,m }_{0,i}=\hat{v}^{\fell-\fe_2,m
			}_{0,i}=\hat{v}^{\fell-\fone, m}_{0,i}$ for $i=1,...,P_{\ell_2}$. 
		}
	}
	
	\For{$n=1:\cN$}{
		\For{$\fell \in \cI$}{
			
			\For{$m=1:M_{\fell}$}{
				
				\If{$\ell_1=0$ and $\ell_2=0$}{
					Compute the EnKF prediction states for $i=1,...,P_{\ell_2}$ 
					$ v^{\fell,m}_{n,i} = \textbf{Prediction}(\hat{v}^{\fell,m}_{n-1,i})$ similar to \eqref{enkf:prediction}. 
					
					Compute the EnKF updated states for $i=1,...,P_{\ell_2}$ $\hat{v}^{\fell,m}_{n,i}=  \textbf{Update}(v^{\fell,m}_{n,i})$ similar to \eqref{enkf:update}.
					
					Compute the EnKF estimator
					$\fDelta \mu_n^{\fell,m}[\varphi]= \sum_{i=1}^{P_{\ell_2}}\frac{\varphi(\vHat_{n, i}^{\fell,m}) }{P_{\ell_2}}.$  
				}
				\ElseIf{$\ell_1>0$ and $\ell_2 = 0 $ }{		
					Compute pairwise coupled prediction states for $i=1,...,P_{\ell_2}$ 
					$ v^{\fell,m}_{n,i}, v^{\fell-\fe_1,m
					}_{n,i}= \textbf{Prediction}(\hat{v}^{\fell,m}_{n-1,i}, \hat{v}^{\fell-\fe_1,m
					}_{n-1,i})$ similar to \eqref{mlenkf:prediction}. 
					
					Compute pairwise coupled updated states for $i=1,...,P_{\ell_2}$ $\hat{v}^{\fell,m}_{n,i},\hat{v}^{\fell-\fe_1,m
					}_{n,i}=  \textbf{Update}(v^{\fell,m}_{n,i}, v^{\fell-\fe_1,m
					}_{n,i})$ similar to \eqref{mlenkf:update}.
					
					Compute the EnKF estimator pairwise coupled in $N_{\ell_1}$
					$\fDelta \mu_n^{\fell,m}[\varphi]= \sum_{i=1}^{P_{\ell_2}}\frac{\varphi(\vHat_{n, i}^{\fell,m}) -\varphi(\vHat_{n, i}^{\fell-\fe_1,m}) }{P_{\ell_2}}.$ }
				\ElseIf{$\ell_1=0$ and $\ell_2>0 $ }{	
					Compute pairwise coupled prediction states for $i=1,...,P_{\ell_2}$ 
					$ v^{\fell,m}_{n,i}, v^{\fell-\fe_2,m
					}_{n,i}= \textbf{Prediction}(\hat{v}^{\fell,m}_{n-1,i}, \hat{v}^{\fell-\fe_2,m
					}_{n-1,i})$ similar to \eqref{mi:prediction1}-\eqref{mi:prediction3}. 
					
					Compute pairwise coupled updated states for $i=1,...,P_{\ell_2}$ $\hat{v}^{\fell,m}_{n,i},\hat{v}^{\fell-\fe_2,m
					}_{n,i}=  \textbf{Update}(v^{\fell,m}_{n,i}, v^{\fell-\fe_2,m
					}_{n,i})$ similar to \eqref{mi:update1}- \eqref{mi:update2}.
					
					Compute the EnKF estimator pairwise coupled in $P_{\ell_2}$
					$\fDelta \mu_n^{\fell,m}[\varphi]= \sum_{i=1}^{P_{\ell_2}}\frac{\varphi(\vHat_{n, i}^{\fell,m}) -\varphi(\vHat_{n, i}^{\fell-\fe_2,m}) }{P_{\ell_2}}.$ }
				\ElseIf{$\ell_1> 0$ and $\ell_2>0 $ }{	
					Compute the four-coupled prediction states for $i=1,...,P_{\ell_2}$ 
					$ v^{\fell,m}_{n,i}, v^{\fell-\fe_1,m}_{n,i}, v^{\fell-\fe_2,m
					}_{n,i}, \hat{v}^{\fell-\fone,m}_{n,i} = \textbf{Prediction}(\hat{v}^{\fell,m}_{n-1,i},\hat{v}^{\fell-\fe_1,m }_{n-1,i},\hat{v}^{\fell-\fe_2,m
					}_{n-1,i},\hat{v}^{\fell-\fone,m}_{n-1,i})$ by \eqref{mi:prediction1}-\eqref{mi:prediction3}. 
					
					Compute the four-coupled updated states for $i=1,...,P_{\ell_2}$ $\hat{v}^{\fell,m}_{n,i},\hat{v}^{\fell-\fe_1,m}_{n,i},\hat{v}^{\fell-\fe_2,m
					}_{n,i},\hat{v}^{\fell-\fone,m }_{n,i} =  \textbf{Update}(v^{\fell,m}_{n,i}, v^{\fell-\fe_1,m}_{n,i}, v^{\fell-\fe_2,m
					}_{n,i}, \hat{v}^{\fell-\fone,m}_{n,i} )$ by \eqref{mi:update1}-\eqref{mi:update2}.
					
					Compute the four-coupled EnKF estimator
					$\fDelta \mu_n^{\fell,m}[\varphi]= \sum_{i=1}^{P_{\ell_2}}\frac{\varphi(\vHat_{n, i}^{\fell,m})-\varphi(\vHat_{n, i}^{\fell-\fe_1,m})  -\varphi(\vHat_{n, i}^{\fell-\fe_2,m}) +\varphi(\vHat_{n, i}^{\fell-\fone,m}) }{P_{\ell_2}}.$ } 
				
			} 
			
		}
		
		Compute the MIEnKF estimator $\mu_n^{MI}[\varphi]=\sum_{ \fell  \in \cI}
		\sum_{m=1}^{M_{\fell}} \frac{\fDelta \mu_n^{\fell,m}[\varphi]}{M_{\fell}}.$ 
	}

\end{algorithm}

\section{MIEnKF complexity}
\label{sec:complexity}
This section presents a cost-versus-accuracy result for the MIEnKF method,
and compares the performance of MIEnKF to MLEnKF and EnKF.

Let us first recall that we restrict ourselves to resolutions of the form 
\[
N_{\ell_1} = N_0\times 2^{\ell_1} \quad \text{and} \quad P_{\ell_2} = P_0\times 2^{\ell_2} \qquad \forall \fell \in \bN_0^2,
\]
for some $N_0, P_0 \in \bN$, and proceed with defining the notion of admissible QoIs:

\begin{definition}[Admissible QoI]
  A Borel-measurable mapping $\varphi: \bR^d \to \bR$ is said to be an admissible QoI if it satisfies the following two integrability conditions for all $n \ge 0$:
  \[
    \bar \mu_n [\varphi] < \infty \quad \text{and} \quad  \mu_n^{\fell}[\varphi] \in L^2(\Omega) \quad \forall \fell \in \bN_0^2.
  \]
\end{definition}  
For any admissible QoI and $\fell \in\bN_0^2$, 
the definition implies that $\fDelta\mu_n^{\fell}[\varphi] \in L^2(\Omega)$,
and we impose the additional assumptions to ensure good performance 
for MIEnKF:

\begin{assumption}\label{ass:Psi2}
  For any admissible QoI $\varphi$ and any $n\ge 0$, the 
  four-coupled EnKF estimator $\fDelta\mu_n^{\fell}[\varphi]$ satisfies the following conditions:
  \[
  \abs{\Ex{\fDelta\mu_n^{\fell}[\varphi]}}  \lesssim N_{\ell_1}^{-1} P_{\ell_2}^{-1},\tag{\textbf{A1}} \label{ass:A1}
  \]
  \[
  \mathbb{V}[\fDelta\mu_n^{\fell}[\varphi]] \lesssim N_{\ell_1}^{-2}P_{\ell_2}^{-2},\tag{\textbf{A2}} \label{ass:A2}
  \]
  and
  \[
  \mathrm{Cost}(\fDelta\mu_n^{\fell}[\varphi]) \eqsim  N_{\ell_1}P_{\ell_2} \tag{\textbf{A3}}. \label{ass:A3} \; \;
  \footnote{The constraint~\eqref{ass:A3} could be stated as a property rather than an assumption, as it holds by construction. Every prediction-update iteration of the coupled EnKF ensembles relating to $\fDelta\mu_n^{\fell}$ costs $\cO(N_{\ell_1} P_{\ell_2})$, cf.~Section~\ref{ssec:mienkfalg}.}
  \]
\end{assumption}

\begin{theorem}[MIEnKF complexity]\label{thm:complexityMIEnKF}
  Let Assumptions~\ref{ass:Psi} and~\ref{ass:Psi2} hold, and for any
  $\epsilon >0$ consider the MIEnKF method with triangular index set 
  $\cI=\{\fell \in \bN_0^2 \mid \ell_1+ \ell_2 \leq L\}$,
  where 
  \begin{equation*}
      L= \max\Big(\lceil \log \epsilon^{-1}+\log\log \epsilon^{-1} \rceil-L_0, \; 1\Big) \quad \text{for some} \quad L_0 \in \bN_0
  \end{equation*}
  and the number of samples 
  \[
      M_{\fell} \eqsim \epsilon^{-2} N_{\ell_1}^{-3/2} P_{\ell_2}^{-3/2} \qquad \fell \in \cI .
  \]

  For any admissible QoI $\varphi$ and $n\ge0$, it then holds that
  \begin{equation}\label{eq:mienkfMSE}
  \Ex{\prt{\mu^{MI}_n[\varphi]-\bar \mu_n [\varphi]}^2}  \lesssim \epsilon^2,
  \end{equation}
  and the computational cost of the MIEnKF estimator satisfies that
  \[
  \mathrm{Cost(\mu_n^{MI}[\varphi])}\eqsim \epsilon^{-2}.
  \]
\end{theorem}
\begin{proof}
  Adding and subtracting $\Ex{\mu^{MI}_n[\varphi]}$ in the mean-squared error, we obtain 
  \[
  \Ex{\prt{\mu^{MI}_n[\varphi] \, \pm \Ex{\mu^{MI}_n[\varphi]} \, -\bar \mu_n [\varphi]}^2}= \Var{\mu^{MI}_n[\varphi]}+\prt{\Ex{\mu^{MI}_n[\varphi]}-\bar \mu_n [\varphi]}^2.
  \]
  For the variance term, the independence of the random variables
  $\{\fDelta \mu_n^{\fell,m}[\varphi]\}_{(\fell, m)}$ and~\eqref{ass:A2} yield
  \[
  \Var{\mu^{MI}_n[\varphi]}= \sum_{ \fell  \in \cI} \sum_{m=1}^{M_{\fell}} \frac{\Var{\fDelta\mu_n^{\fell, m}[\varphi]}}{M_{\fell}} \lesssim  \sum_{ \fell  \in \cI} M_{\fell}^{-1} N_{\ell_1}^{-2}P_{\ell_2}^{-2} \lesssim \epsilon^2.
  \]
  For the squared bias term,~\eqref{ass:A1} and the multi-index telescoping properties of the MIEnKF estimator imply that
  \[
  \prt{\Ex{\mu^{MI}_n[\varphi]}-\bar \mu_n [\varphi]}^2 \leq \prt{\sum_{ \fell  \notin \cI} \Ex{\fDelta\mu_n^{\fell}[\varphi]}}^2 \lesssim  \prt{\sum_{ \fell  \notin \cI} N_{\ell_1}^{-1} P_{\ell_2}^{-1}}^2.
  \]
  The mean-squared error bound~\eqref{eq:mienkfMSE} follows by
  \begin{equation*}
    \begin{split}
      \sum_{ \fell  \notin \cI} N_{\ell_1}^{-1} P_{\ell_2}^{-1} \lesssim
      \sum_{\ell_1 + \ell_2 > L}^{\infty}2^{-(\ell_1+\ell_2)} = \sum_{k=L+1}^\infty (k+1)2^{-k} \eqsim 2^{-L}L \eqsim \epsilon
    \end{split}
  \end{equation*}
  and 
  \[
  \mathrm{Cost}(\mu_n^{MI}[\varphi]) = \sum_{ \fell  \in \cI}
  M_{\fell} \mathrm{Cost} (\fDelta \mu_n^{\fell} [\varphi])  \eqsim \sum_{ \fell  \in \cI}
  M_{\fell} N_{\ell_1}P_{\ell_2}  \eqsim \epsilon^{-2}. \; \; \footnote{Optimal choices for $L$ and $M_{\fell}$ in Theorem~\ref{thm:complexityMIEnKF}
    can be obtained by solving a constrained optimization problem, analogously as elaborated on for the MLMC method in the seminal paper~\cite{giles2008}.}
  \]
\end{proof}

\begin{remark}
  For comparison, we briefly recall the cost-versus-accuracy results
  for EnKF and MLEnKF. For any $\epsilon>0$ and
  sufficient regularity
  \[
  \begin{split}
    &\norm{(\mu_n^{N,P}-\bar \mu_n)[\varphi]}_p \lesssim \epsilon, \qquad \mathrm{(EnKF)}\\
    &\norm{(\mu_n^{ML}-\bar \mu_n)[\varphi]}_p \lesssim \epsilon, \qquad \mathrm{(MLEnKF)}
  \end{split}
  \]
  with the computational cost bounded by
  \[
  \begin{split}
    &\mathrm{Cost}(\mu_n^{N,P}[\varphi])\eqsim \epsilon^{-3},\\
    &\mathrm{Cost}(\mu_n^{ML}[\varphi])\eqsim \epsilon^{-2} \abs{\log(\epsilon)}^3.
  \end{split}
  \]
  For more details, see~\cite{hoel2020multilevel}.
\end{remark}

\begin{remark}
  An alternative to~\eqref{ass:A2} that is more aligned
  with assumption made for the existing convergence results for MLEnKF
is to assume that 
  \[ 
  \norm{\fDelta\mu_n^{\fell}[\varphi]}_p \lesssim N_{\ell_1}^{-1}P_{\ell_2}^{-1} \tag{\textbf{A2$^\star$}}, \label{ass:A2star}
  \]
  for $p\ge 2$. Then, using the same aforementioned index set $\cI$ and $L$ and with a slight change in the
  sample size $M_{\fell} \eqsim \epsilon^{-2} N_{\ell_1}^{-4/3}P_{\ell_2}^{-4/3}$, the MIEnKF estimator
  satisfies 
  \begin{equation}\label{mienkf:LpError}
    \norm{(\mu_n^{MI}-\bar \mu_n)[\varphi]}_p \lesssim \epsilon,
  \end{equation}
  with the asymptotic MIEnKF cost bounded by $\mathcal{O}(\epsilon^{-2})$.
  This can be proved similarly as the case of Theorem~\ref{thm:complexityMIEnKF}, where the
  $L_p$-norm of the statistical error can be bounded using the Marcinkiewicz-Zygmund inequality:
  \begin{equation*}
    \begin{split}
      \norm{\prt{\mu^{MI}_n-\bar \mu_n}[\varphi]}_p &\leq \norm{\mu^{MI}_n[\varphi]-\E[\mu^{MI}_n][\varphi]}_p+\norm{\E[\mu^{MI}_n][\varphi]-\bar \mu_n[\varphi]}_p\\
      & \lesssim  \sum_{ \fell  \in \cI} M_{\fell}^{-1/2} \norm{\fDelta\mu_n^{\fell, m}[\varphi] }_p+ \sum_{ \fell  \notin \cI} \abs{\E[\fDelta\mu_n^{\fell,m}[\varphi]]}.
    \end{split}
  \end{equation*}
  We note that \eqref{ass:A2} is a weaker assumption than \eqref{ass:A2star} since
  \begin{equation}\label{ass:A2p=2}
    \Var{\fDelta\mu_n^{\fell}[\varphi]} \leq \norm{\fDelta\mu_n^{\fell}[\varphi]}_2^2 \le \norm{\fDelta\mu_n^{\fell}[\varphi]}_p^2.
  \end{equation}
\end{remark}

\begin{remark}
  Under more general settings, Assumption~\ref{ass:Psi2} may be transformed into
  \[ \abs{\E[\fDelta\mu_n^{\fell}[\varphi]] }  \lesssim N_{\ell_1}^{-\alpha_1} P_{\ell_2}^{-\alpha_2},\tag{$\mathcal{A}_1$} \label{ass:A1gen}\]
  \[\norm{\fDelta\mu_n^{\fell}[\varphi]}_p \lesssim N_{\ell_1}^{-\beta_1}P_{\ell_2}^{-\beta_2},\tag{$\mathcal{A}_2$} \label{ass:A2gen}\]
  \[ \mathrm{Cost}(\fDelta\mu_n^{\fell}[\varphi]) \eqsim
    N_{\ell_1}^{\gamma_1}P_{\ell_2}^{\gamma_2}
    \tag{$\mathcal{A}_3$}, \label{ass:A3gen}\] for some
  $\alpha_1, \alpha_2, \beta_1, \beta_2, \gamma_1, \gamma_2 >0$.  The
  construction of an efficient MIEnKF estimator may then lead to a
  differently shaped (possibly even non-triangular) index set $\cI$, a
  different sequence of number of samples
  $\{M_{\fell}\}_{\fell \in \cI}$, and other common ratios for the
  geometric sequences $\{N_{\ell_1}\}$ and $\{P_{\ell_2}\}$. The
  problem of optimizing the set $\cI$ may be recast as a knapsack
  problem, which is a well-studied optimization problem with many
  available solution algorithms, cf.~\cite{abdo2016}
  and~\cite[equation~(21)]{abdo2018}.  For instance, the
  approach developed for approximations of multi-index
  Monte Carlo applied to McKean-Vlasov dynamics in~\cite{abdo2018}
  defines the set by
  \[
    \cI=\{(\ell_1,\ell_2)\in \bN_0^2: (\alpha_1+\gamma_1)\ell_1+(\alpha_2+\gamma_2)\ell_2 \leq L\}.
  \]
\end{remark}

\section{Numerical examples}
\label{sec:numerics}
This section presents a numerical comparison of MIEnKF with the
EnKF and MLEnKF methods outlined in Section~\ref{sec:problem}. Three
problems will be considered: the Ornstein-Uhlenbeck (OU)
process, a stochastic differential equation (SDE) with a double-well (DW) potential, 
and Langevin dynamics~\cite{ hoel2016, hoel2020multilevel, ballesio2020, apte2007sampling, christensen2012forecasting}.

We consider SDE on the general form
\begin{equation}\label{sde:genform}
  du=-U'(u)dt+\sigma dW_t,   
\end{equation}
with a constant diffusion coefficient $\sigma =0.5$
and two types of potential functions:

$(i)\; \:U(u)=u^2/2, \qquad \qquad  \qquad  \qquad   \mbox{  \textbf{(OU)}}$

$(ii)\: U(u)=u^2/4 +1/(4u^2+2), \qquad \mbox{\textbf{(DW)}}$.\\

The numerical discretizations of~\eqref{sde:genform} are computed using
the Milstein numerical scheme\footnote{Note that in all three examples considered, the SDEs are with constant diffusion terms. For such SDEs, the Milstein scheme coincides with the Euler-Maruyama scheme.} with uniform timestep $\Delta t = 1/N$
for any $N\ge 1$.  The observations of the process $u$
are equally spaced with 
observation time interval $\tau=1$, observation operator $H=1$,  
$\Gamma=0.1$ and the QoI $\varphi(x)=x$.

To numerically verify assumptions~\eqref{ass:A1}
and~\eqref{ass:A2star}, the following rates are estimated from
$S$ independent copies of $\fDelta \mu_{n}^{\fell}[\varphi]$:
\begin{equation*}
  \begin{split}
    \abs{\E[\fDelta\mu_n^{\fell} [\varphi]]} &\approx \Big| \sum_{i=1}^{S} \frac{\fDelta \mu_{n,i}^{\fell}[\varphi]}{S} \Big|,\\
    \norm{\fDelta\mu_{n}^{\fell}[\varphi]}_2& \approx\sqrt{\frac{1}{S}\sum_{i=1}^{S} \abs{\fDelta \mu_{n,i}^{\fell}[\varphi]}^2 }.
  \end{split}
\end{equation*}
We analyze the convergence rates of the methods by computing the
time-averaged root-mean-squared error (RMSE).
\[
\mbox{RMSE}:=\sqrt{ \frac{1}{S(\mathcal{N}+1)}\sum_{i=1}^{S}\sum_{n=0}^{\mathcal{N}}\abs{\mu^{*}_{n,i} [\varphi] - \bar{\mu}_n [\varphi]}^2},
\]
where $\{\mu_{.,i}^{*}[\varphi]\}_{i=1}^{S}$ are independent copies of
$\mu_{.}^{*}[\varphi]$ for the specific methods (EnKF, MLEnKF, and
MIEnKF).

\subsection{Reference solutions and computer architecture}
Since dynamics $\Psi$ is linear for the OU problem, the reference solution $\bar{\mu}_n[\varphi]$ can be computed exactly
using the Kalman filter. However, the reference solution for the DW problem, which involves
nonlinear dynamics, must be approximated. This solution
is computed using the deterministic mean-field EnKF
algorithm, cf.~\cite[Appendix C]{hoel2020multilevel}. A
pseudoreference solution for the final test problem based on Langevin
dynamics is computed by the sample average of
$S=180$ independent simulations of the
MIEnKF estimator at the tolerance
$\epsilon=2^{-11}$ using the following parameters:
\begin{equation*}\label{par:refmienkf}
  \begin{split}
    L&= \lceil L_*+\log_2(L_*)\rceil -1, \qquad \text{with} \quad L_*= \lceil \log_2(\epsilon^{-1})\rceil -1, \\
    N_{\ell_1} &= 4 \times 2^{\ell_1}, \\
    P_{\ell_2} &= 30 \times 2^{\ell_2}, \\
    M_{\fell} &=
    \begin{cases}
      6 \times \lceil \epsilon^{-2} N_{\ell_1}^{-3/2}P_{\ell_2}^{-3/2} \rceil  & \text{if} \quad \ell_1 =0 \mbox{ and } \ell_2 =0,\\
      90 \times \lceil \epsilon^{-2} N_{\ell_1}^{-3/2}P_{\ell_2}^{-3/2} \rceil & \text{if} \quad  1 \le \ell_1+\ell_2 \le  L.
    \end{cases}
  \end{split}
\end{equation*}

The numerical simulations were computed in parallel on 18 cores on
an Intel(R) Xeon(R) CPU E5-2680 v2 20-core processor with 128 GB RAM.
The computer code was written in the Julia programming
language~\cite{Julia-2017}, and it can be downloaded from
\url{https://github.com/GaukharSH/mienkf}.

\subsection{Ornstein-Uhlenbeck process} 
\label{ssec:ou}
\begin{figure}[h!]
	\centering
	\hspace*{-1cm}
	\includegraphics[width=1.2\textwidth]{{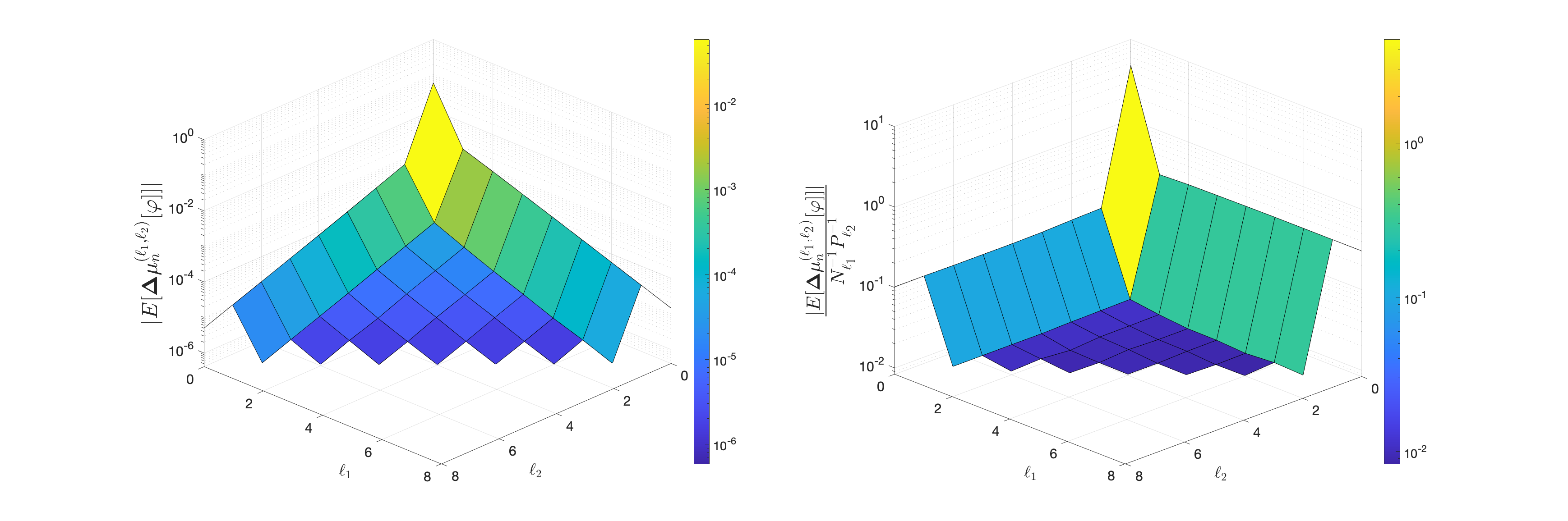}}
	\hspace*{-1cm}
	\includegraphics[width=1.2\textwidth]{{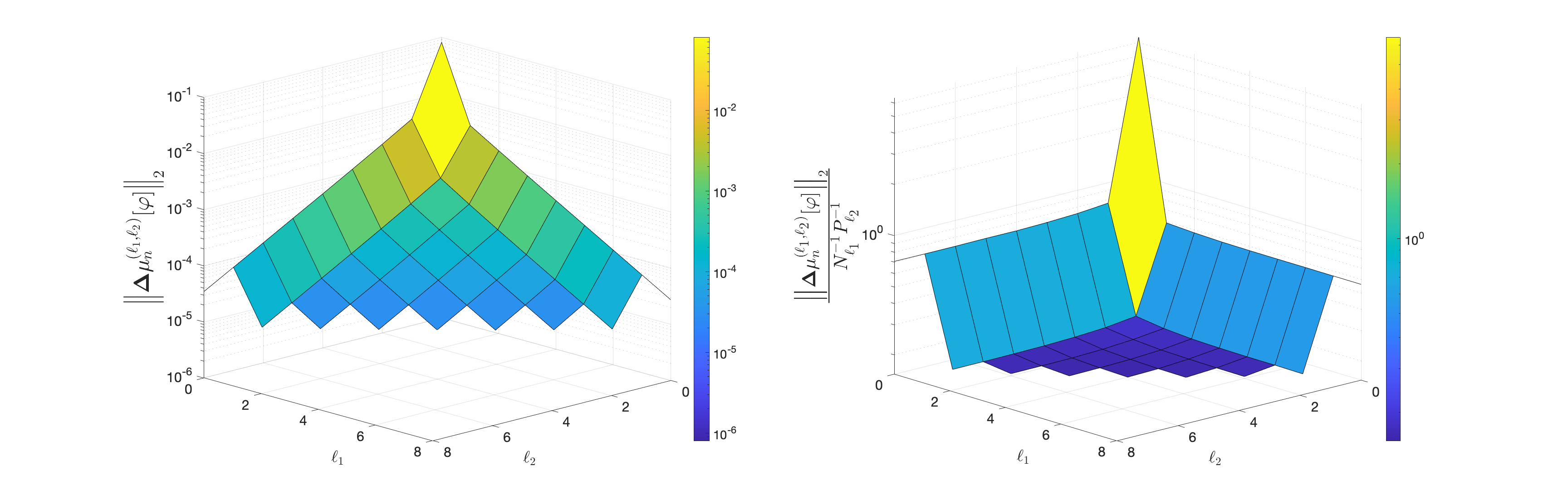}}
	\caption{\textbf{Ornstein-Uhlenbeck problem. Estimates based on $S=10^6$ independent
			runs (Section~\ref{ssec:ou}).} Top row: Numerical evidence of
		assumption~\eqref{ass:A1} over $\cN=20$ observation times when using
		$N_{\ell_1}=4\times 2^{\ell_1}$ and $P_{\ell_2}=20 \times
		2^{\ell_2}$.  Bottom row: Similar plots for the verification of assumption~\eqref{ass:A2star}.}
	\label{fig:OUconjRates}
\end{figure}

We consider the SDE~\eqref{sde:genform} with the (OU) potential
function and initial condition $u(0)\sim N(0, \Gamma)$. Convergence rates~\eqref{ass:A1}
and~\eqref{ass:A2star} shown in Figure~\ref{fig:OUconjRates}  were estimated by the Monte Carlo
method using $S=10^6$ independent samples of
$\fDelta\mu_n^{\fell} [\varphi]$.  In the figure,  the left panel shows the weak and
$L_2$ convergence rates over $\cN=10$ observation times with
$(\ell_1+ \ell_2) \in [0,7]$, and the right panel shows the ratio of
the rates to $N_{\ell_1}^{-1}P_{\ell_2}^{-1}$. The plane-like flatness
of the right panel for $(\ell_1+\ell_2) \in [1,7]$ validates the
said rate assumptions.
    
When conducting runtime-versus-accuracy convergence tests for an input
tolerance $\epsilon>0$, we set the parameters of the
respective methods as follows:
\begin{equation}\label{par:enkf}
  \textbf{EnKF:} \quad  P = \lceil15 \epsilon^{-2}\rceil \quad \text{and} \quad  N  = \lceil \epsilon^{-1} \rceil,
\end{equation}
\begin{equation}\label{par:mlenkf}
  \textbf{MLEnKF:} \quad \left\{
  \begin{split}
    L&= \lceil \log_2(\epsilon^{-1})\rceil-1, \\
    N_\ell &= 2 \times 2^{\ell}, \\
    P_\ell &= 10 \times 2^\ell, \\
    M_\ell &=
    \begin{cases}
      2 \times \lceil \epsilon^{-2} L^2 2^{-3}\rceil & \text{if} \quad \ell =0,\\
      \lceil \epsilon^{-2} L^2 2^{-2\ell-3}\rceil & \text{if} \quad  1 \le \ell \le  L,
    \end{cases}
  \end{split}
  \right.
\end{equation}
and  
\begin{equation}\label{par:mienkf}
  \textbf{MIEnKF:} \quad \left\{
  \begin{split}
    L&= \lceil L_*+\log_2(L_*)\rceil -1, \quad
    \text{with} \quad L_*= \lceil \log_2(\epsilon^{-1})\rceil -1,\\
    N_{\ell_1} &= 4 \times 2^{\ell_1}, \\
    P_{\ell_2} &= 30 \times 2^{\ell_2}, \\
    M_{\fell} &=
    \begin{cases}
      6 \times \lceil \epsilon^{-2} N_{\ell_1}^{-3/2}P_{\ell_2}^{-3/2} \rceil & \text{if} \quad \ell_1 =0 \mbox{ and } \ell_2 =0,\\
      120 \times \lceil \epsilon^{-2} N_{\ell_1}^{-3/2}P_{\ell_2}^{-3/2} \rceil  & \text{if} \quad  1 \le \ell_1+\ell_2 \le  L.
    \end{cases}
  \end{split}
  \right.
\end{equation}

\begin{figure}[h!]
  \centering
  \includegraphics[width=0.47\textwidth]{{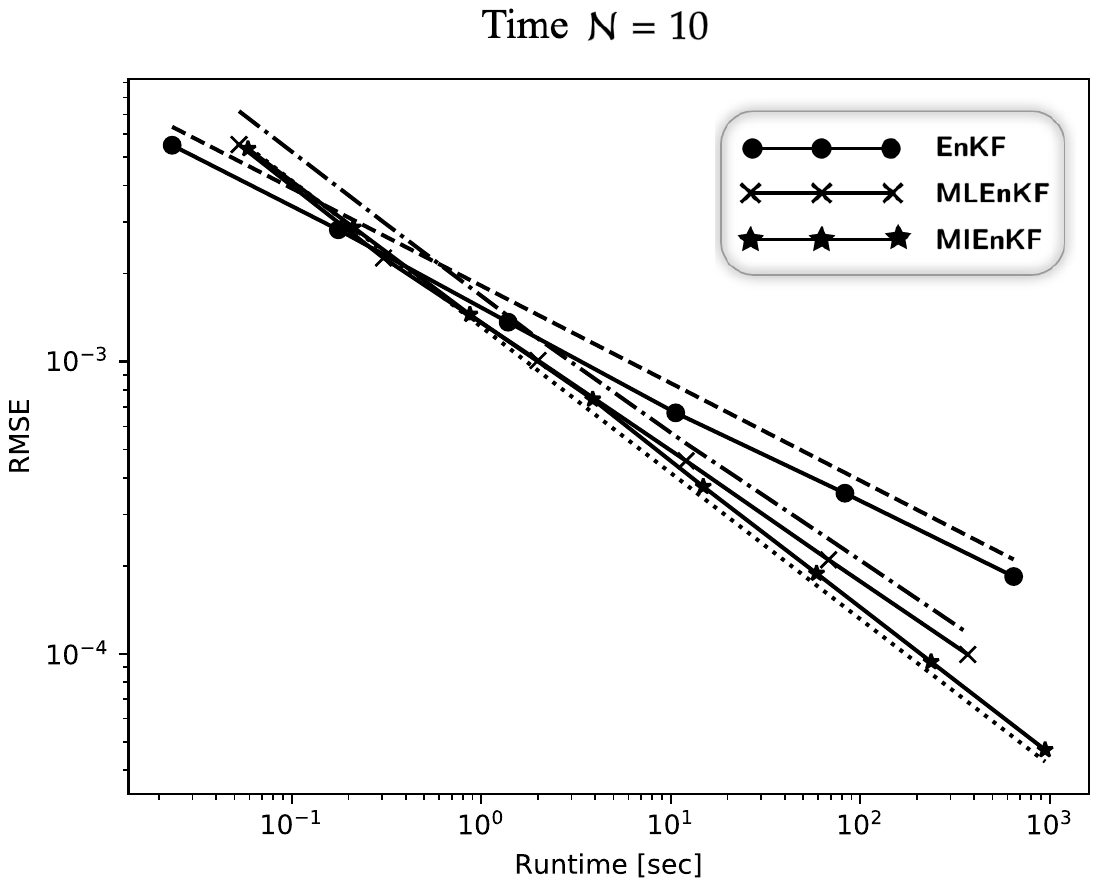}}
  \includegraphics[width=0.47\textwidth]{{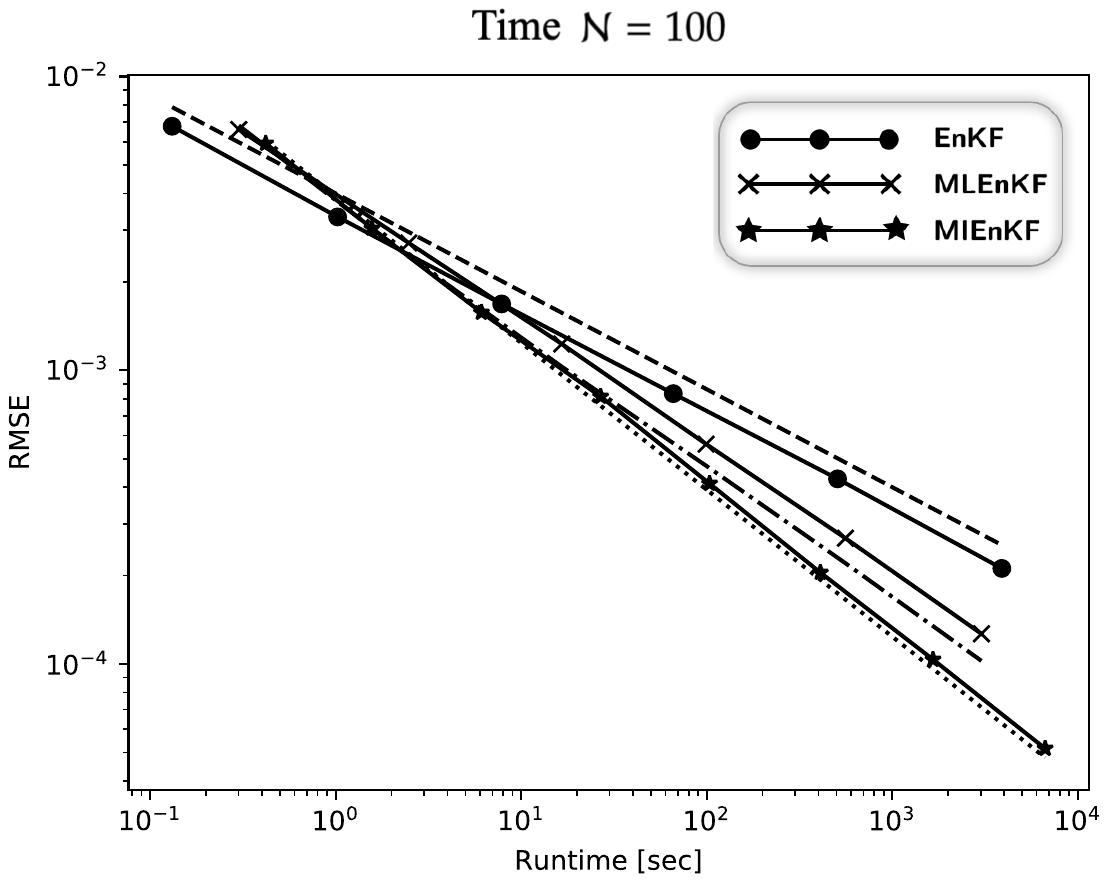}}
  \vspace*{-2cm}
  \caption{\textbf{Ornstein-Uhlenbeck problem. Estimates based on $S=100$ independent
      runs (Section~\ref{ssec:ou}).}  Comparison of the runtime versus root-mean-squared error (RMSE)
    for mean over observation times $\cN=10$ (left) and $\cN=100$ (right). The solid-crossed line represents MLEnKF and the
    dot-dashed line is a fitted
    $\cO(\log(10+\mathrm{Runtime})^{1/3}\mathrm{Runtime}^{-1/2})$
    reference line. The solid-asterisk line represents the MIEnKF and the
    dotted line is a fitted $\cO(\mathrm{Runtime}^{-1/2})$ reference
    line. The solid-bulleted line represents EnKF and the dashed line
    is a fitted $\cO(\mathrm{Runtime}^{-1/3})$ reference line.}
  \label{fig:OUconvRates}
\end{figure}

For a sequence of predefined tolerances $\epsilon = [2^{-4}, 2^{-5},
  2^{-6}, 2^{-7}, 2^{-8}, 2^{-9}]$ for EnKF and MLEnKF, and  $\epsilon = [2^{-4}, 2^{-5},
  2^{-6}, 2^{-7}, 2^{-8}, 2^{-9}, 2^{-10}, 2^{-11}]$ for MIEnKF,
Figure~\ref{fig:OUconvRates} shows the runtime against the RMSE for
the three methods over observation times of 
$\mathcal{N}=10$ and $\mathcal{N}=100$ estimated using $S=100$
independent runs. MIEnKF outperforms EnKF and MLEnKF
for sufficiently small tolerances, and the
complexity rate agrees with the theory. 


\subsection{Double-well SDE} 
\label{ssec:dw} 
\begin{figure}[h!]
	\centering
	\hspace*{-1cm}
	\includegraphics[width=1.2\textwidth]{{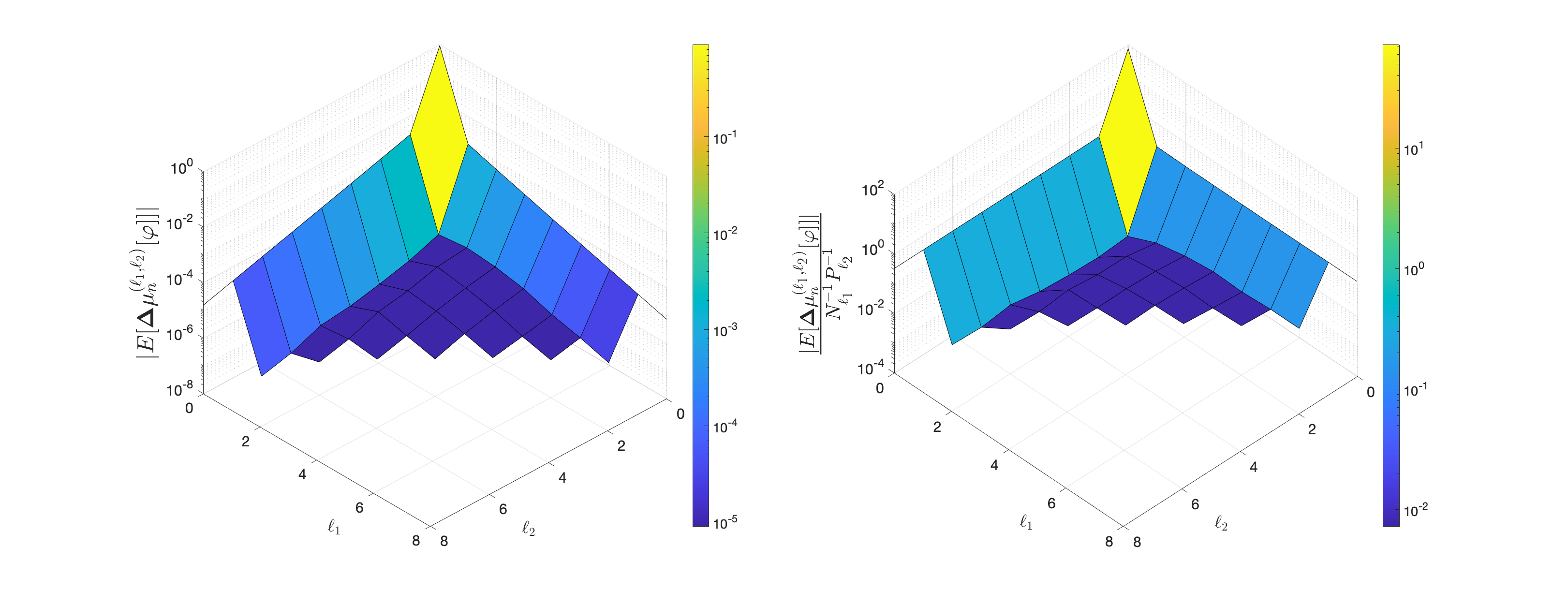}}
	\hspace*{-1cm}
	\includegraphics[width=1.2\textwidth]{{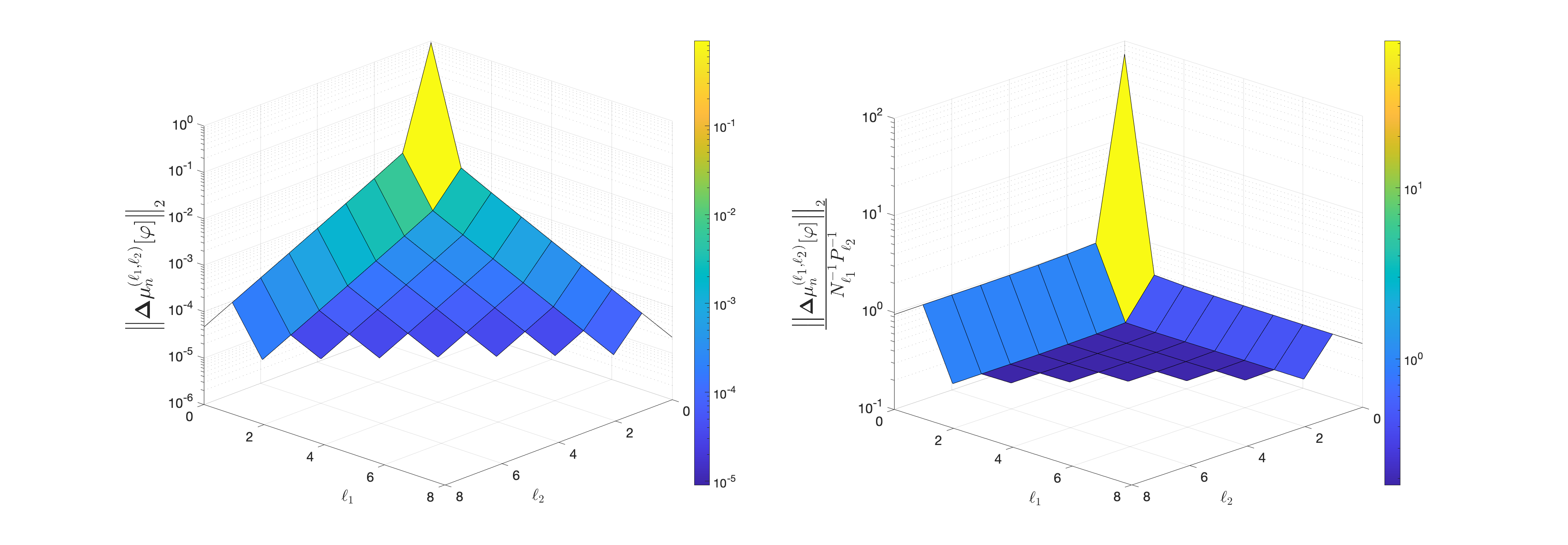}}
	\caption{\textbf{Double Well problem. Estimates based on $S=10^6$ independent runs (Section~\ref{ssec:dw}).}
		Top row: Numerical evidence of assumption~\eqref{ass:A1}
		for $\cN=10$ observation times when using
		$N_{\ell_1}=4\times 2^{\ell_1}$ and
		$P_{\ell_2}=20 \times 2^{\ell_2}$.  Bottom row: Similar plots for
		verifying assumption~\eqref{ass:A2star}.  }
	\label{fig:DWconjRates}
\end{figure}

\begin{figure}[h!]
  \centering
  \includegraphics[width=0.47\textwidth]{{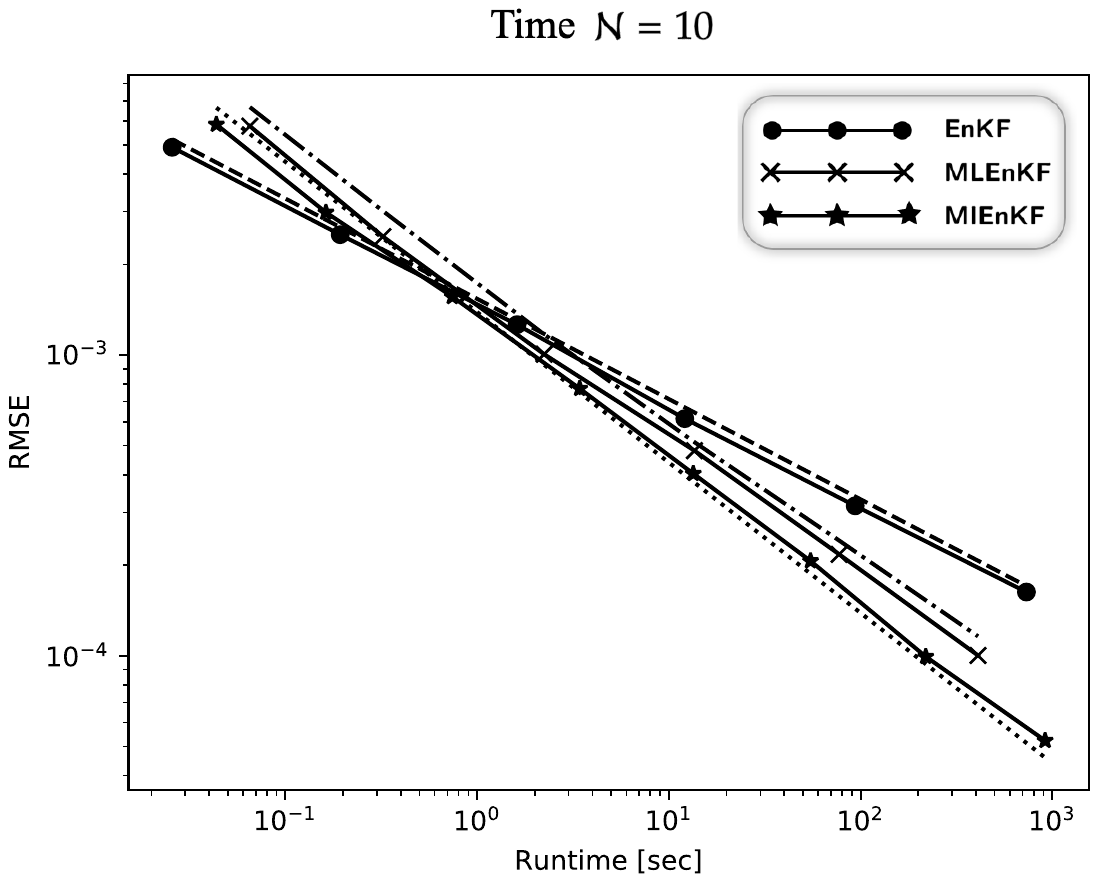}}
  \includegraphics[width=0.47\textwidth]{{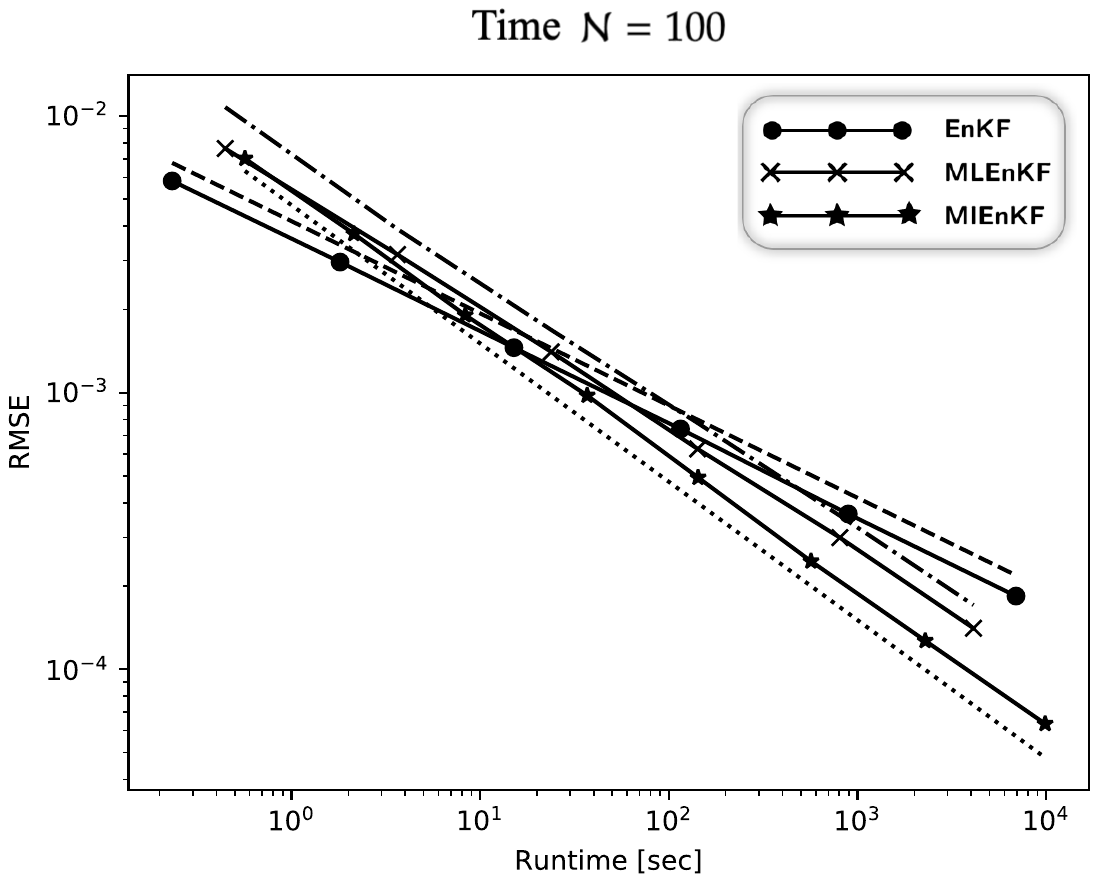}}
 \vspace*{-2cm}
  \caption{\textbf{Double Well problem. Estimates based on $S=100$ independent runs (Section~\ref{ssec:dw}).}
    Similar plots as those shown in Figure~\ref{fig:OUconvRates}.
  }
  \label{fig:DWconvRates}
\end{figure}

We consider the SDE~\eqref{sde:genform} with the DW potential
function and $u(0)\sim N(0, \Gamma)$.  Similar to the OU case,
Figure~\ref{fig:DWconjRates} provides numerical evidence of the
conjecture rates under assumptions~\eqref{ass:A1} and
\eqref{ass:A2star}. For the same predefined $\epsilon-$inputs with the
same degrees of freedom setting as in the example of OU, the performance of
the three methods were compared in terms of runtime against RMSE for
observation times $\mathcal{N}=10$ and $\mathcal{N}=100$ and
estimated over $S=100$ independent runs
(Figure~\ref{fig:DWconvRates}). We observe that MIEnKF outperforms EnKF and
MLEnKF for small RMSE.


\subsection{Langevin SDE}\label{ssec:langevin}
In the last example, we consider the two-dimensional stochastic
Langevin dynamics
\begin{equation}\label{sde:Langevin}
  \begin{split}
    dX_t&=V_tdt, \\
    dV_t&=-U'(X_t)dt - \kappa V_t dt + (2\kappa T)^{1/2} dW_t,
  \end{split}
\end{equation}
where $X_t$ and $V_t$ denotes the particle position and velocity,
respectively, $U(X)$ is the previously introduced DW potential,
$\kappa=2^{-5}\times\pi^2$ is the viscosity and $T=1$ is the temperature.
To improve the pairwise coupling between particles, we used the first-order
symplectic Euler splitting scheme~\cite{muller2015improving}.
The initial conditions are provided by
$X_0 \sim N(0, \Gamma)$ and $V_0 \sim N(0, \Gamma)$ with
$X_0$ and $V_0$ being independent. Further, based on initial test runs,
the method parameters are set to
\begin{equation*}\label{par: enkf}
  \textbf{EnKF:} \quad  P = \lceil 10 \epsilon^{-2}\rceil \quad \text{and} \quad  N  = \lceil \epsilon^{-1} \rceil,
\end{equation*}
\begin{equation*}\label{par:mlenkf}
    \textbf{MLEnKF:} \quad \left\{ \begin{split}
        L&= \lceil \log_2(\epsilon^{-1})\rceil-1,\\
        N_\ell &= 2 \times 2^{\ell},\\
        P_\ell &= 8 \times 2^\ell, \;\\
    M_\ell &=
    \begin{cases}
      2 \times \lceil \epsilon^{-2} L^2 2^{-2}\rceil & \text{if} \quad \ell =0,\\
      \lceil \epsilon^{-2} L^2 2^{-2\ell-2}\rceil & \text{if} \quad  1 \le \ell \le  L,
    \end{cases}
  \end{split}\right.
\end{equation*}
and 
\begin{equation*}\label{par:mienkf}
  \textbf{MIEnKF:} \quad \left\{
    \begin{split}
        L&= \lceil L_*+\log_2(L_*)\rceil -1,\; \quad
        \text{with} \quad L_*= \lceil \log_2(\epsilon^{-1})\rceil -1 \\
        N_{\ell_1} &= 4 \times 2^{\ell_1},\\
        P_{\ell_2} &= 20 \times 2^{\ell_2},\\
    M_{\fell} &=
    \begin{cases}
      6 \times \lceil \epsilon^{-2} N_{\ell_1}^{-3/2}P_{\ell_2}^{-3/2} \rceil & \text{if} \quad \ell_1 =0 \mbox{ and } \ell_2 =0,\\
      50 \times \lceil \epsilon^{-2} N_{\ell_1}^{-3/2}P_{\ell_2}^{-3/2} \rceil  & \text{if} \quad  1 \le \ell_1+\ell_2 \le  L.
    \end{cases}
  \end{split}\right.
\end{equation*}

\begin{figure}[h!]
  \centering
  \includegraphics[width=0.47\textwidth]{{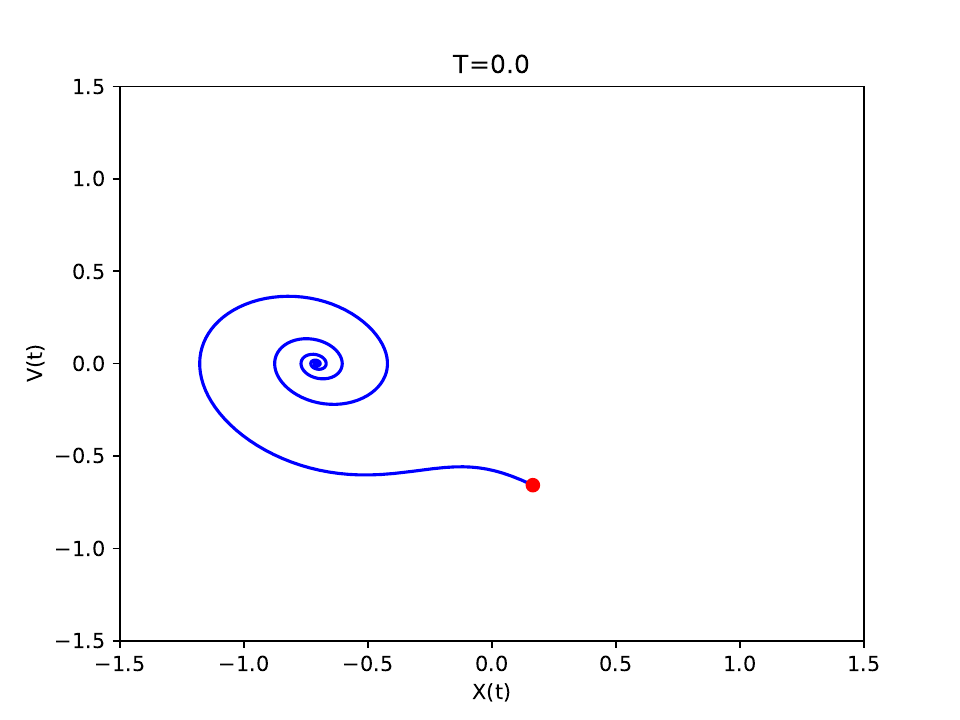}}
  \includegraphics[width=0.47\textwidth]{{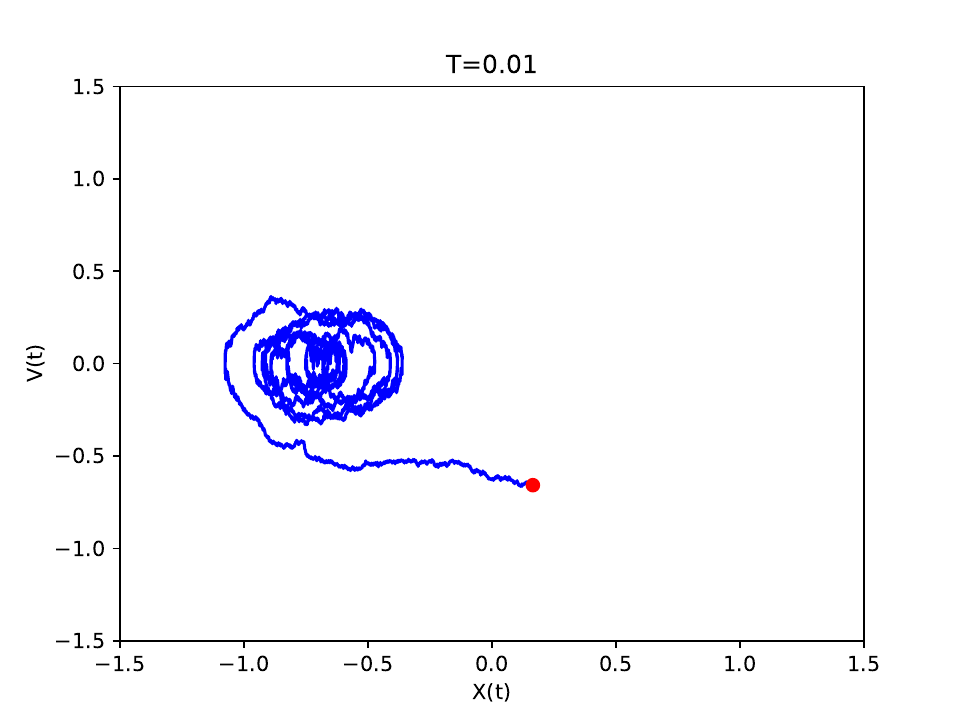}}\\
  \includegraphics[width=0.47\textwidth]{{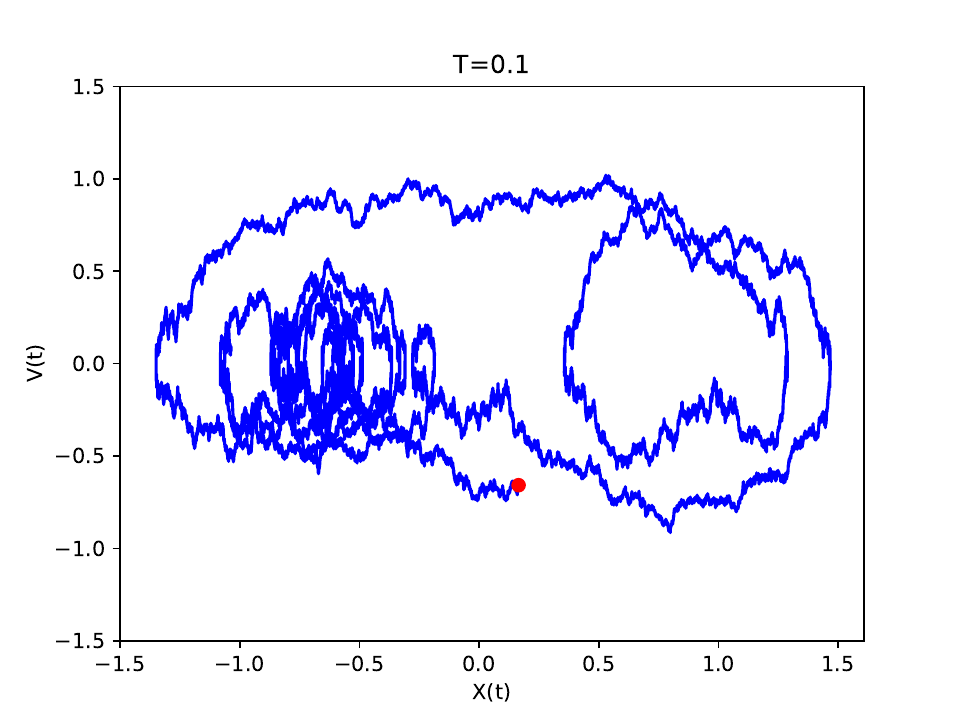}}
  \includegraphics[width=0.47\textwidth]{{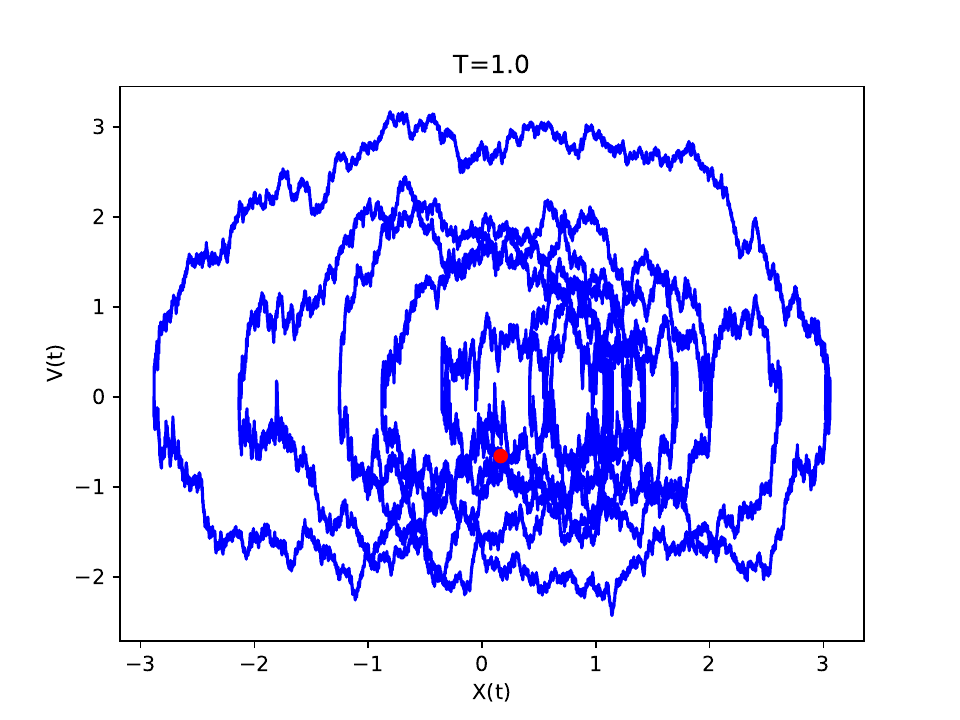}}\\
  \caption{\textbf{Time evolution of a solution to Langevin dynamics
      with different temperature $T$ values}. The symplectic Euler scheme is used
    up to final time $\cN=50$. The red dot represents the initial
    value.}
  \label{fig:LangevinEvolutionT50Split}
\end{figure}

\begin{figure}[h!]
  \centering
  \includegraphics[width=0.47\textwidth]{{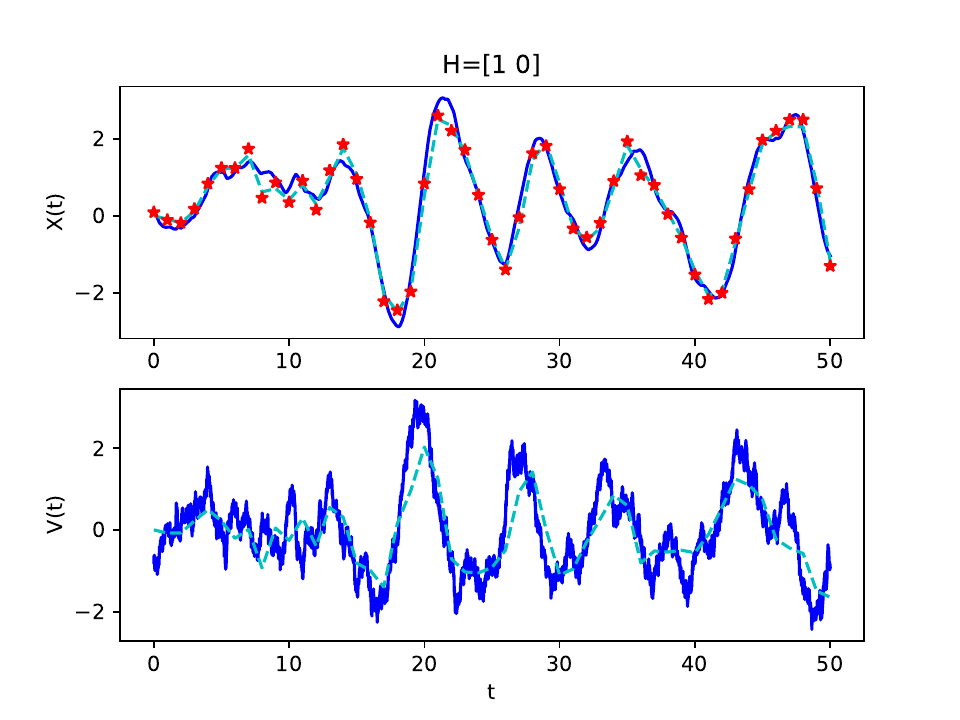}}
  \includegraphics[width=0.47\textwidth]{{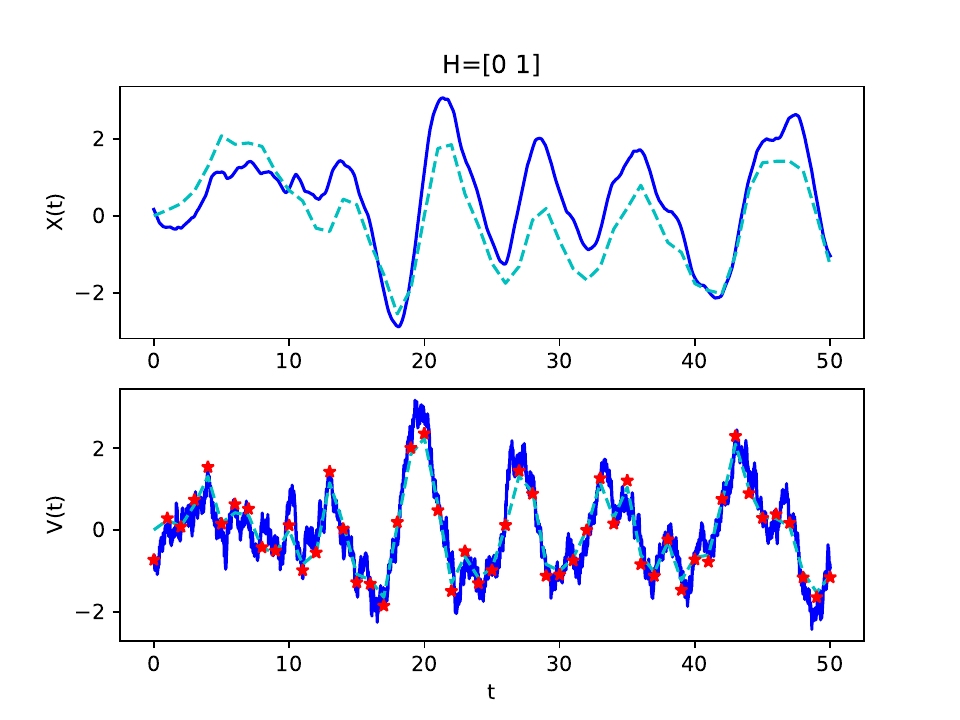}}
   \includegraphics[width=0.47\textwidth]{{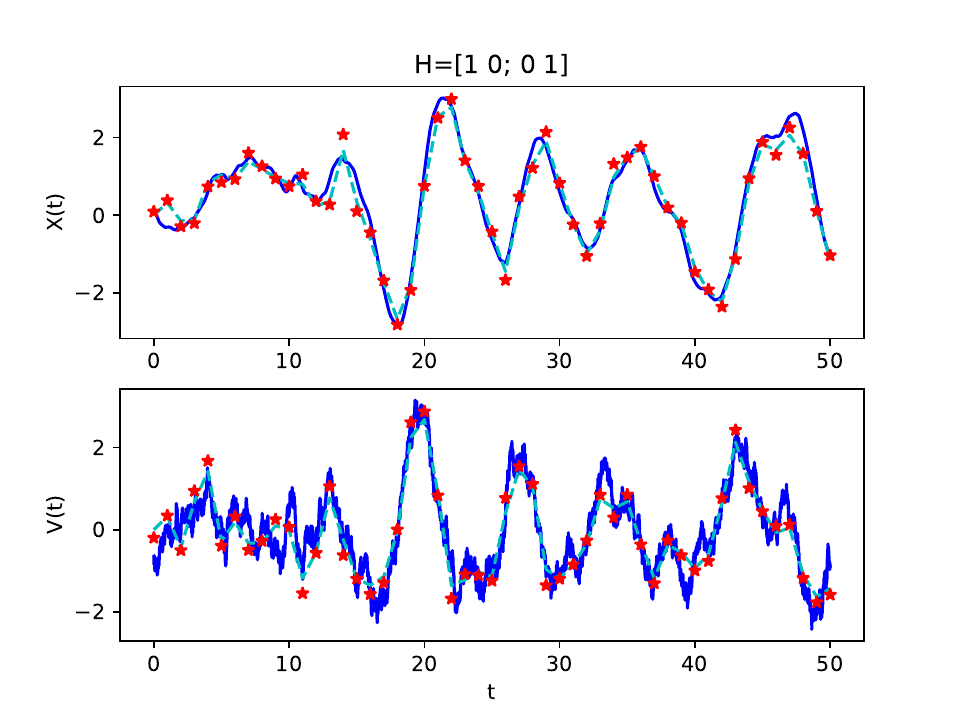}}
  \caption{\textbf{Data assimilation for Langevin dynamics
      given different observation operators $H$,
      Section~\ref{ssec:langevin}}. The blue solid line represents the
    truth, and the red stars are the observations. The cyan dashed lines
    represent MIEnKF mean. The final time $\cN=50$ is used with the
    observation timestep $\tau=1.0$}
  \label{fig:LangevinMeanTruthObsT50SplitH11}
\end{figure}

\begin{figure}[h!]
	\centering
	\begin{subfigure}[t]{4.5in}
		\centering
		\hspace*{-1cm}
		\includegraphics[width=1.2\textwidth]{{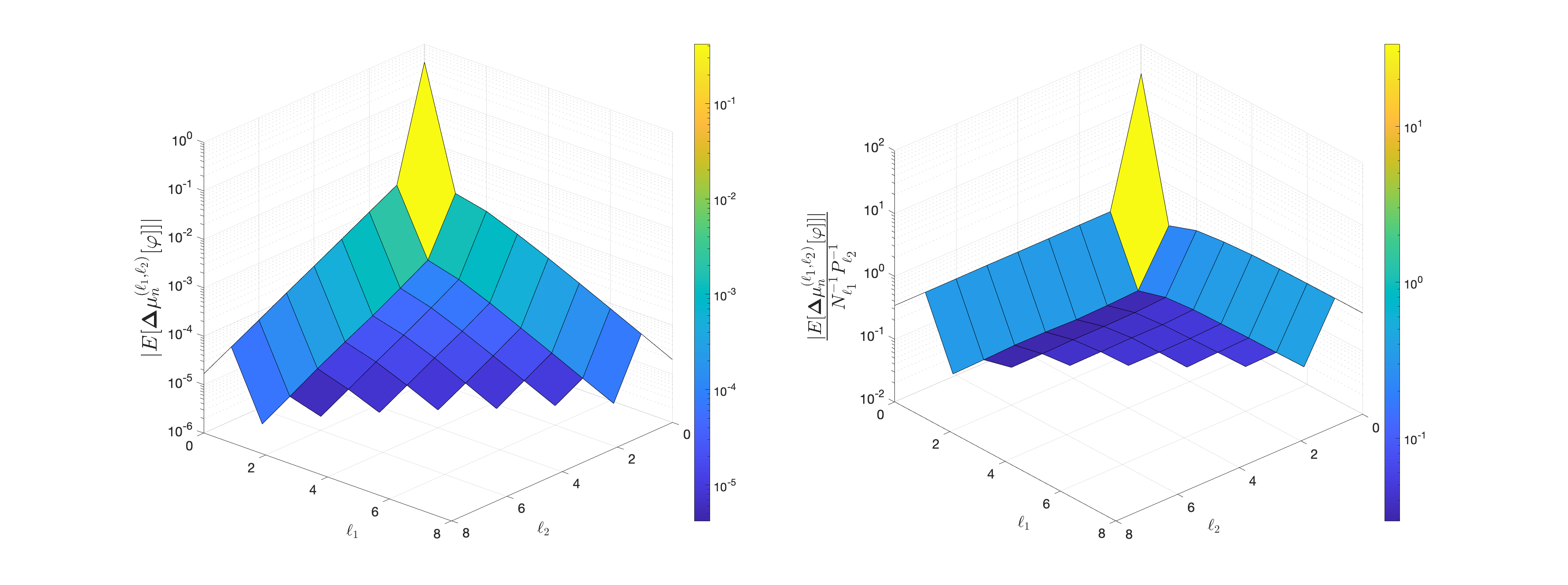}}
		\hspace*{-1cm}
		\includegraphics[width=1.2\textwidth]{{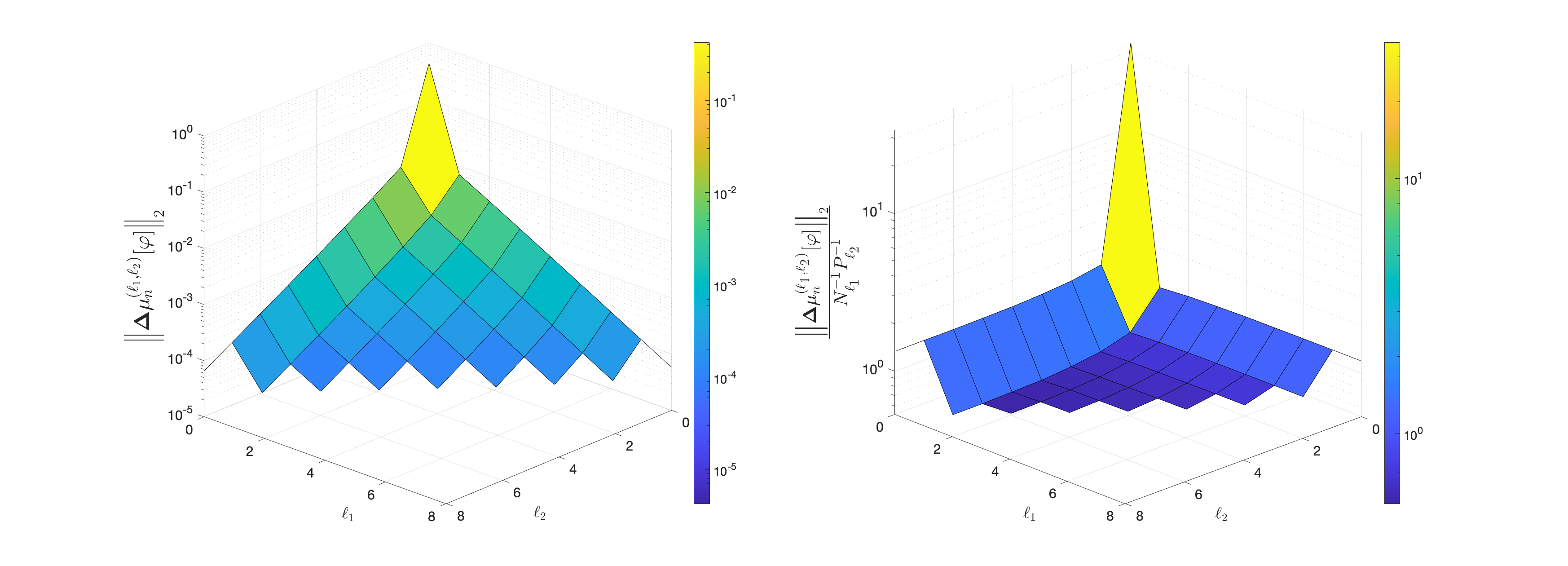}}
		\caption{The particle position $X_t$}\label{fig:LangevinConjRates(a)}
	\end{subfigure}
	\\
	\begin{subfigure}[t]{4.5in}
		\centering
		\hspace*{-1cm}
		\includegraphics[width=1.2\textwidth]{{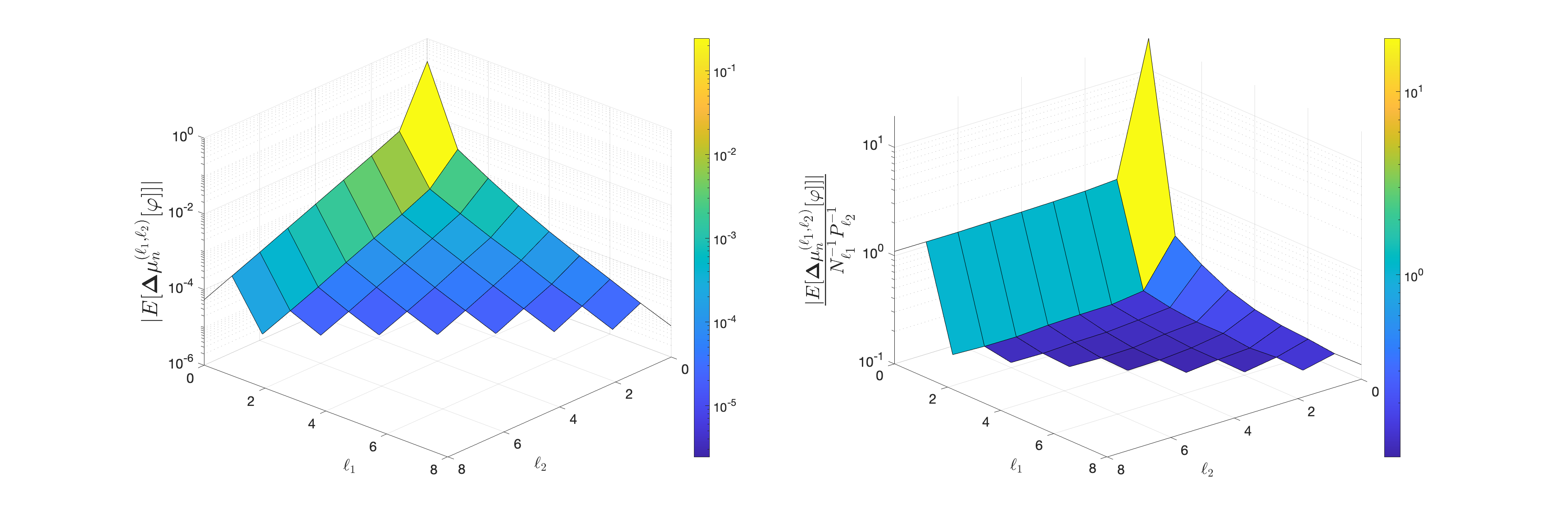}}
		\hspace*{-1cm}
		\includegraphics[width=1.2\textwidth]{{figures/mSXlanH10par1e6Cbar.png}}
		\caption{The particle velocity $V_t$}\label{fig:LangevinConjRates(b)}
	\end{subfigure}
	\caption{\textbf{Langevin dynamics with partial observations, $H=[1 \; 0]$. Estimates based on $S=10^6$
			independent runs  (Section~\ref{ssec:langevin}).} Top row in each subfigure: Numerical
		evidence of assumption~\eqref{ass:A1} for $\cN=10$ observation
		times when using $N_{\ell_1}=4\times 2^{\ell_1}$ and
		$P_{\ell_2}=20 \times 2^{\ell_2}$.  Bottom row in each subfigure:
		Similar plots for verifying
		assumption~\eqref{ass:A2star}.}\label{fig:LangevinConjRatesH10}
\end{figure}

\begin{figure}[h!]
	\centering
	\begin{subfigure}[t]{4.5in}
		\centering
		\hspace*{-1cm}
		\includegraphics[width=1.2\textwidth]{{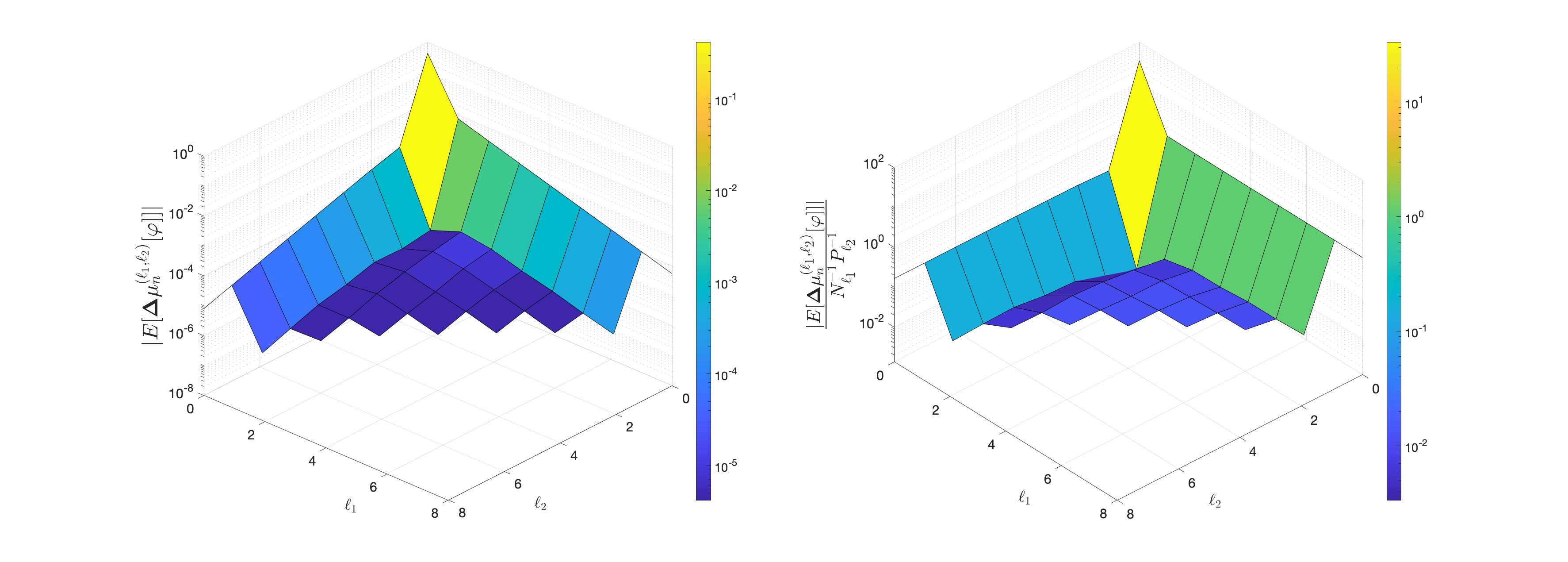}}
		\hspace*{-1cm}
		\includegraphics[width=1.2\textwidth]{{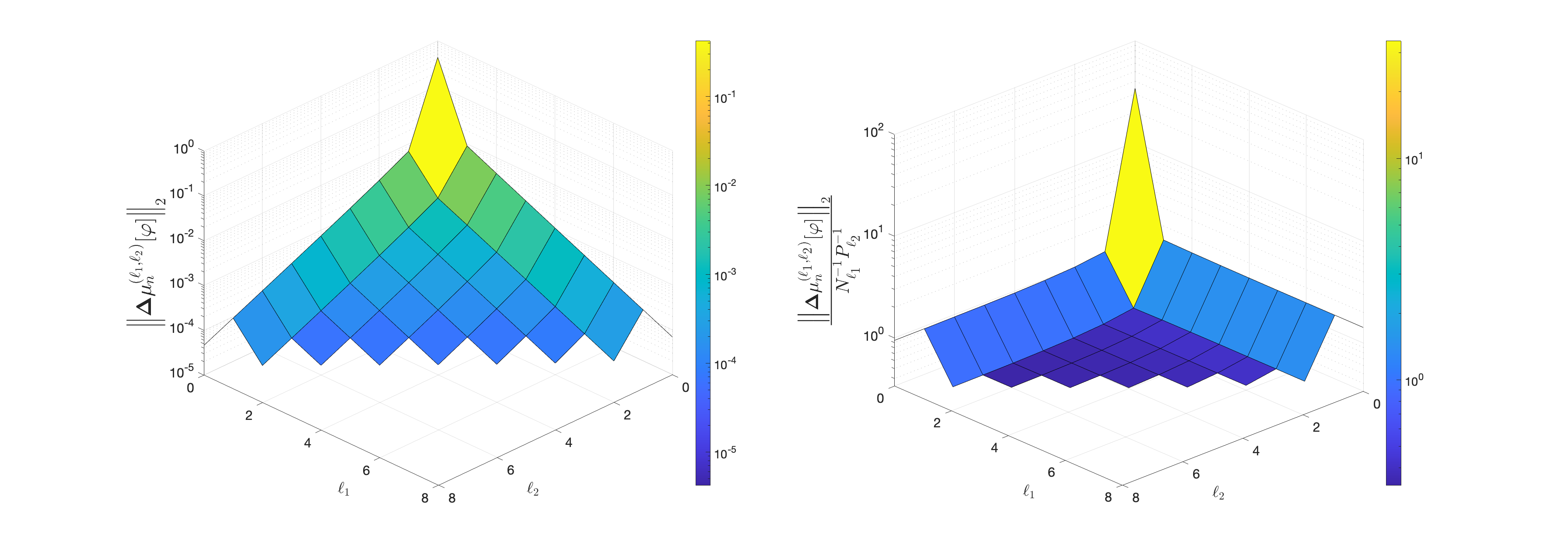}}
		\caption{The particle position $X_t$}\label{fig:LangevinConjRates(a)}
	\end{subfigure}
	\\
	\begin{subfigure}[t]{4.5in}
		\centering
		\hspace*{-1cm}
		\includegraphics[width=1.2\textwidth]{{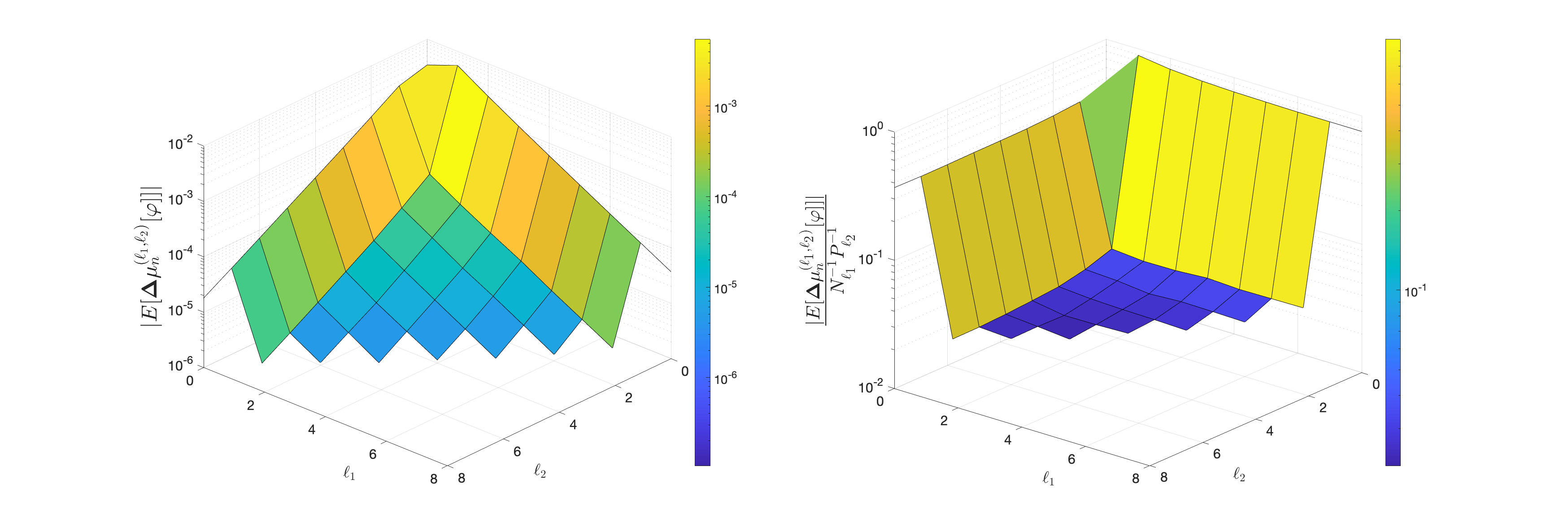}}
		\hspace*{-1cm}
		\includegraphics[width=1.2\textwidth]{{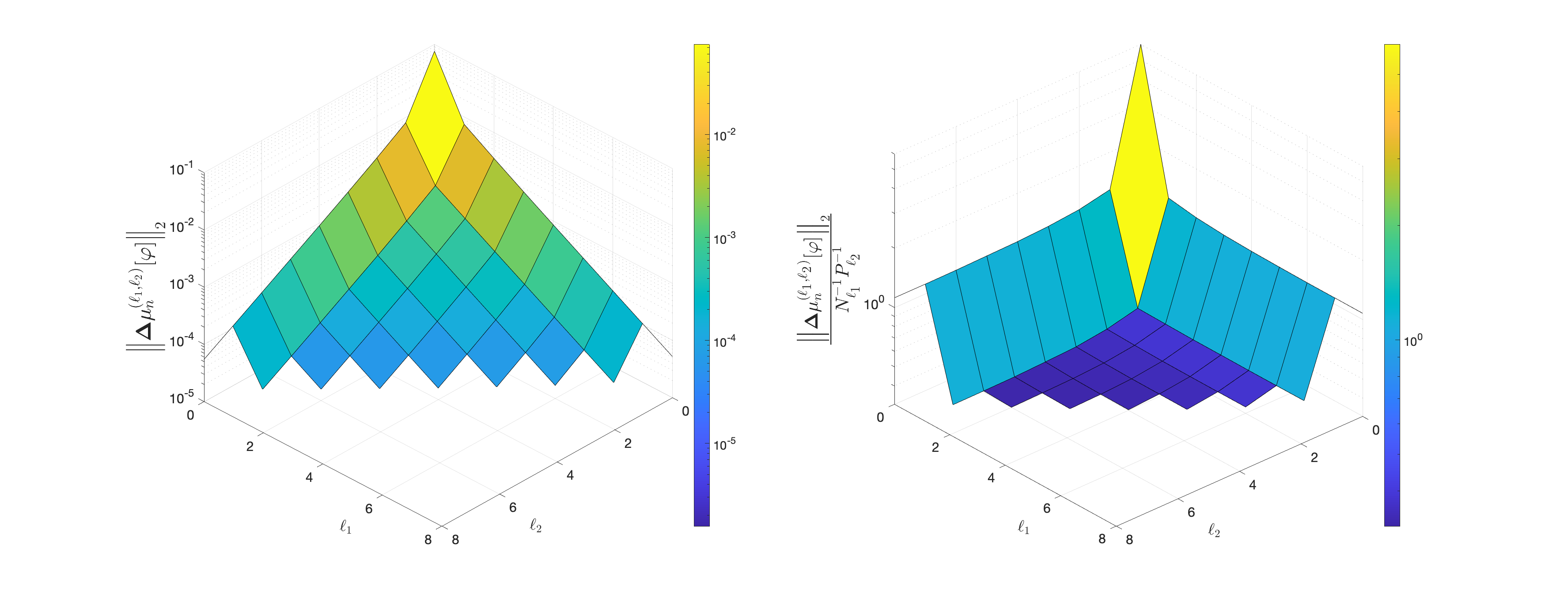}}
		\caption{The particle velocity $P_t$}\label{fig:LangevinConjRates(b)}
	\end{subfigure}
	\caption{\textbf{Langevin dynamics with full observations, $H=[1\;0; 0\; 1]$. Estimates based on $S=10^6$
			independent runs (Section~\ref{ssec:langevin}).} Similar plots as those shown in
		Figure~\ref{fig:LangevinConjRatesH10}.}
	\label{fig:LangevinConjRatesH01}
\end{figure}

\begin{figure}[h!]
	\centering
	\includegraphics[width=0.47\textwidth]{{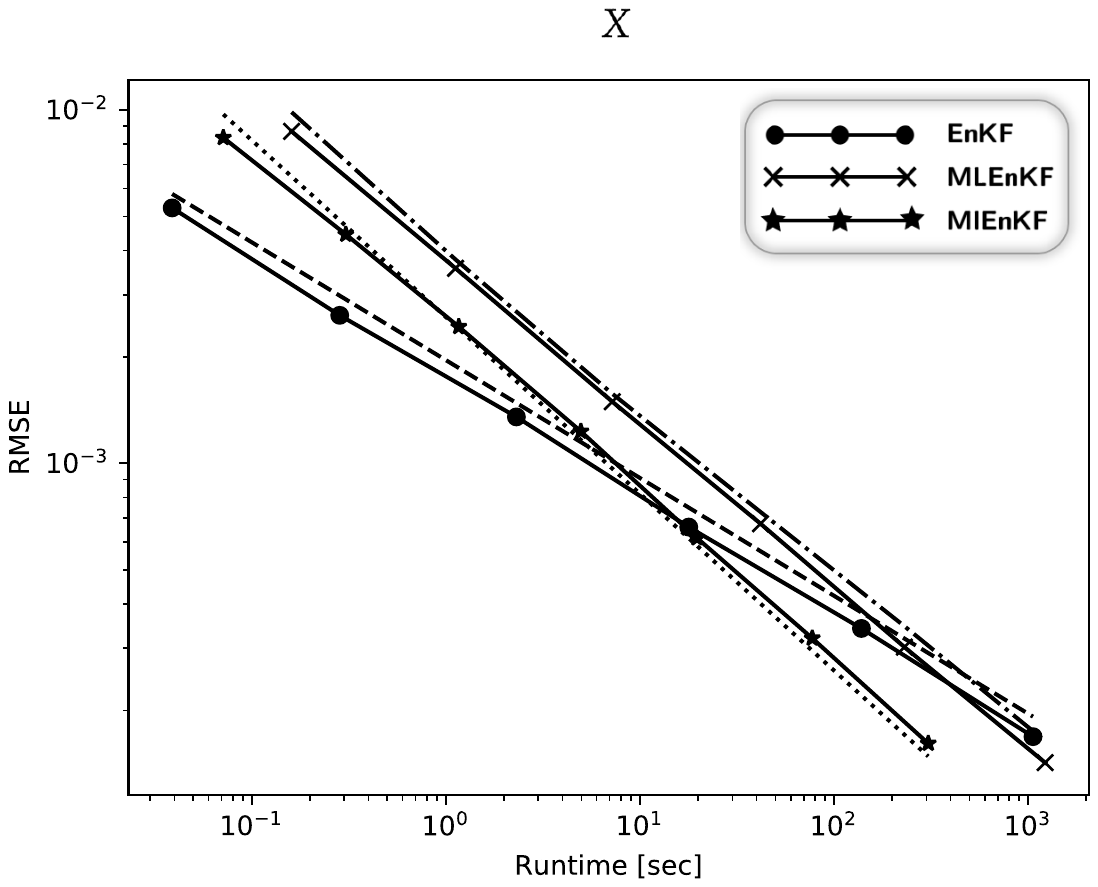}}
	\includegraphics[width=0.47\textwidth]{{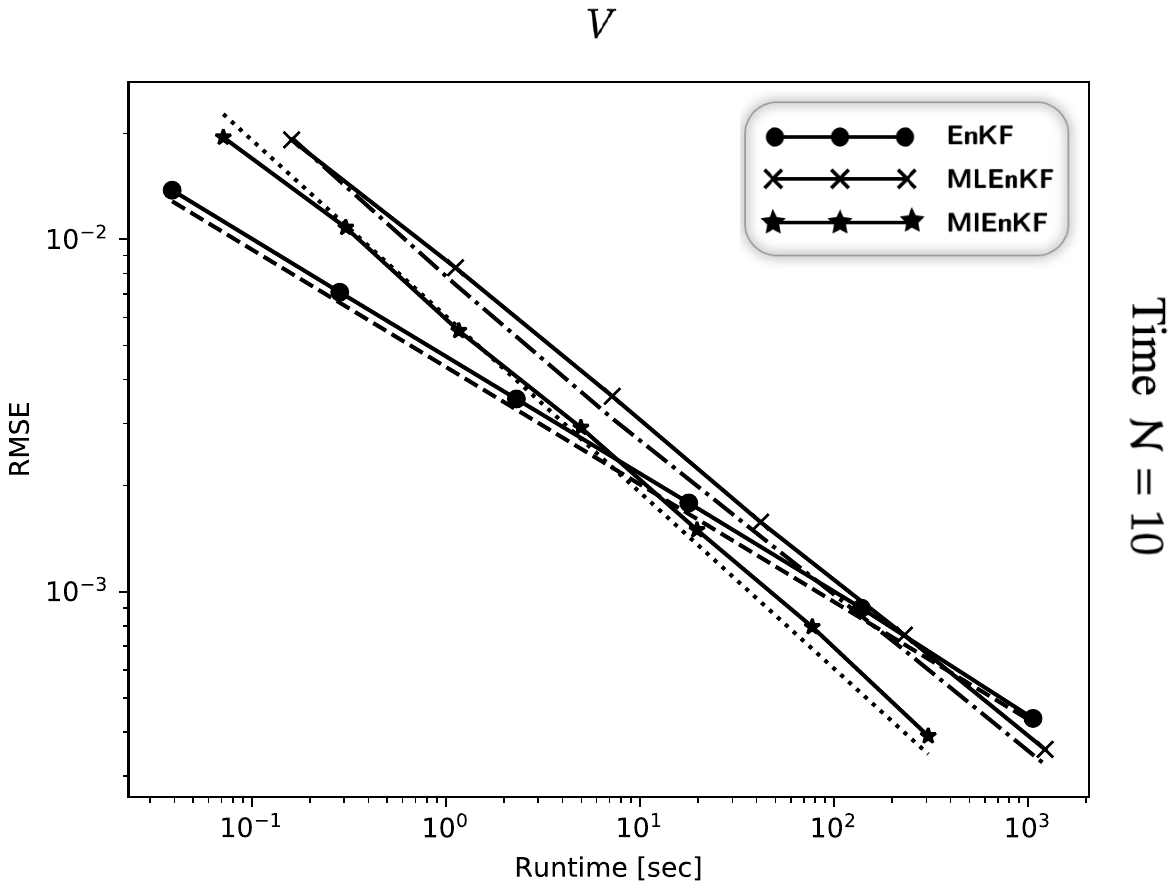}}\\
	\vspace*{-4cm}
	\includegraphics[width=0.47\textwidth]{{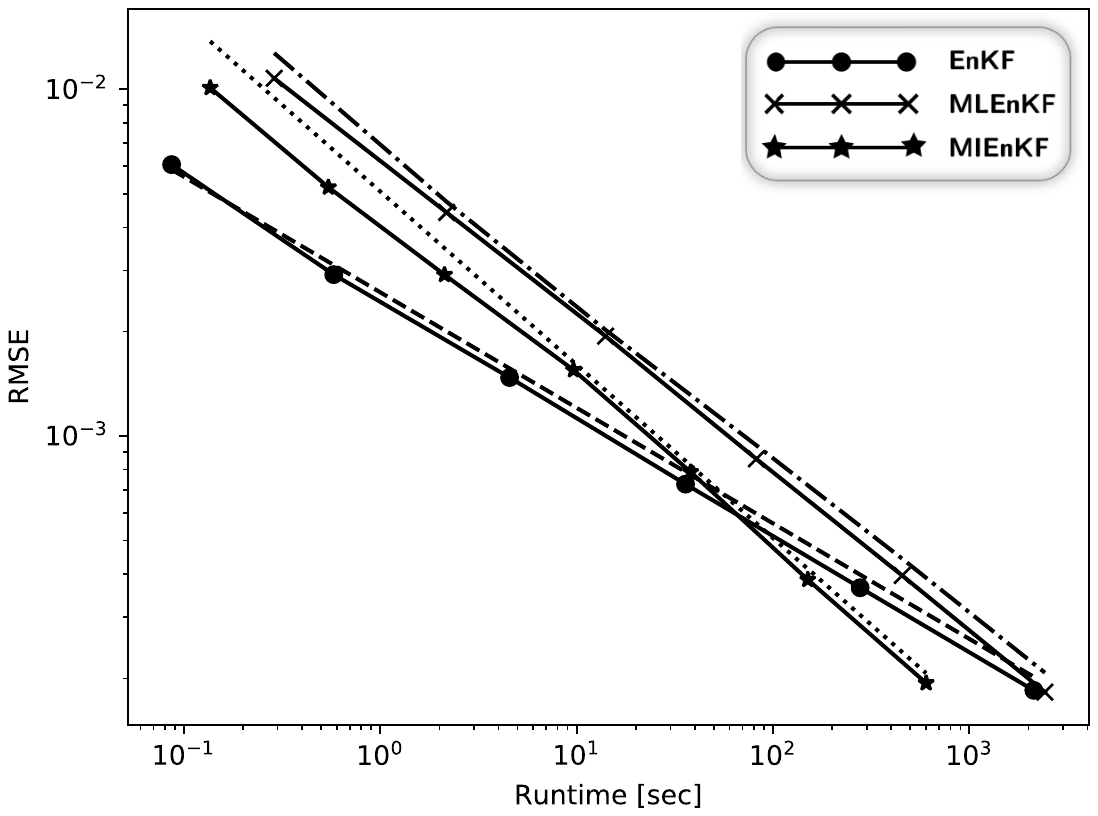}}
	\includegraphics[width=0.47\textwidth]{{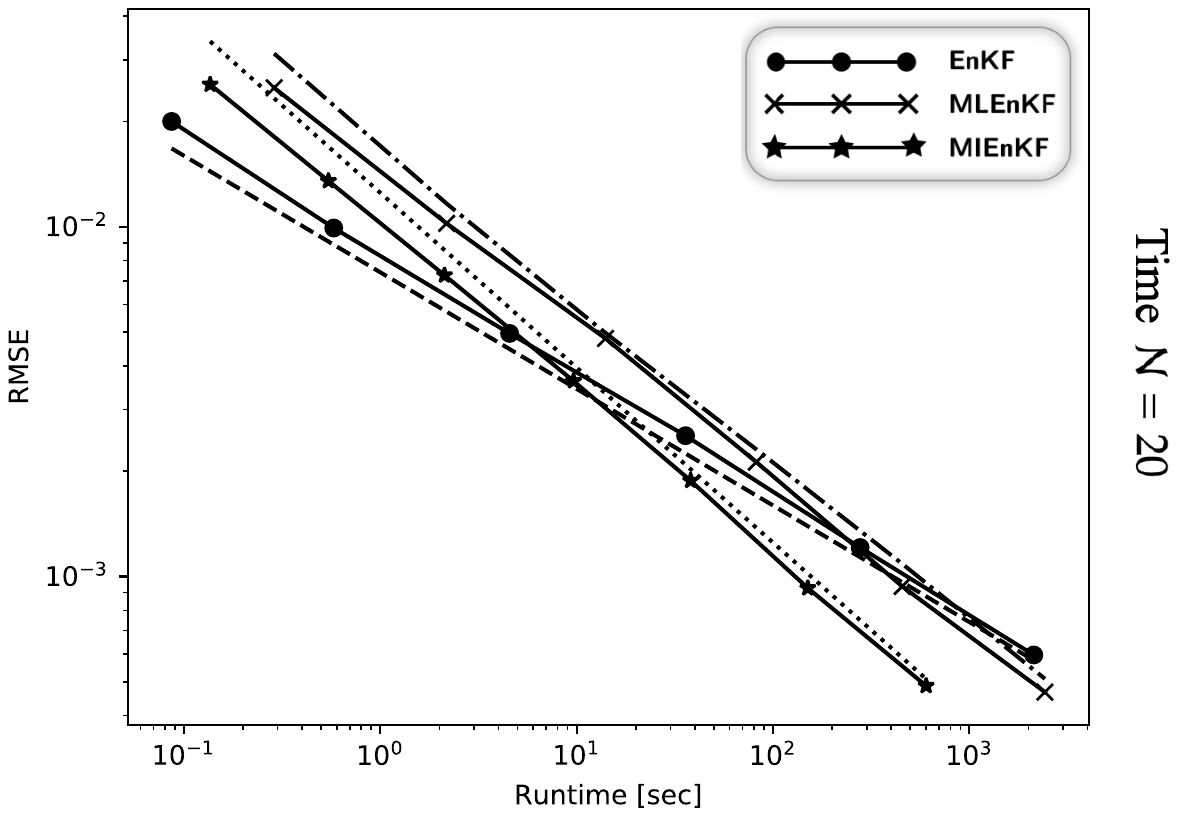}}
	\vspace*{-2cm}
	\caption{\textbf{Langevin dynamics with partial observations, $H=[1 \; 0]$. Estimates based on $S=90$
			independent runs (Section~\ref{ssec:langevin}).}
		\textit{Top row}: Comparison of the runtime versus RMSE for the
	    mean of the component $X$ (left) and the component $V$ (right) over $\cN=10$ observation
		times. The
		solid-crossed line represents MLEnKF and the dot-dashed line is a
		fitted
		$\cO(\log(10+\mathrm{Runtime})^{1/3}\mathrm{Runtime}^{-1/2})$
		reference line.  The solid-asterisk line represents the MIEnKF and the
		dotted line is a fitted $\cO(\mathrm{Runtime}^{-1/2})$ reference
		line. The solid-bulleted line represents EnKF and the dashed line
		is a fitted $\cO(\mathrm{Runtime}^{-1/3})$ reference line.
		\textit{Bottom row}: Similar plots for $\cN=20$ observation
		times.}
	\label{fig:loglogRatesLangevinH10sepT10}
\end{figure}

\begin{figure}[h!]
	\centering
	\includegraphics[width=0.47\textwidth]{{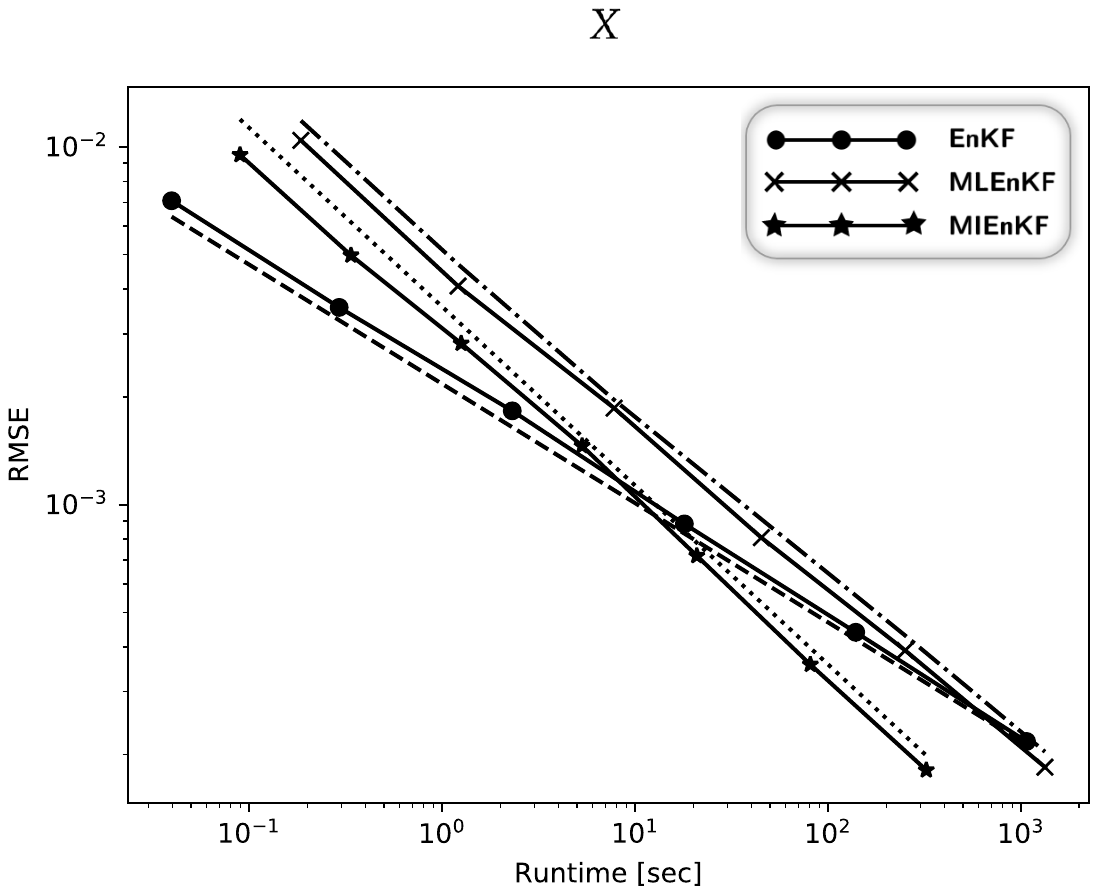}}
	\includegraphics[width=0.47\textwidth]{{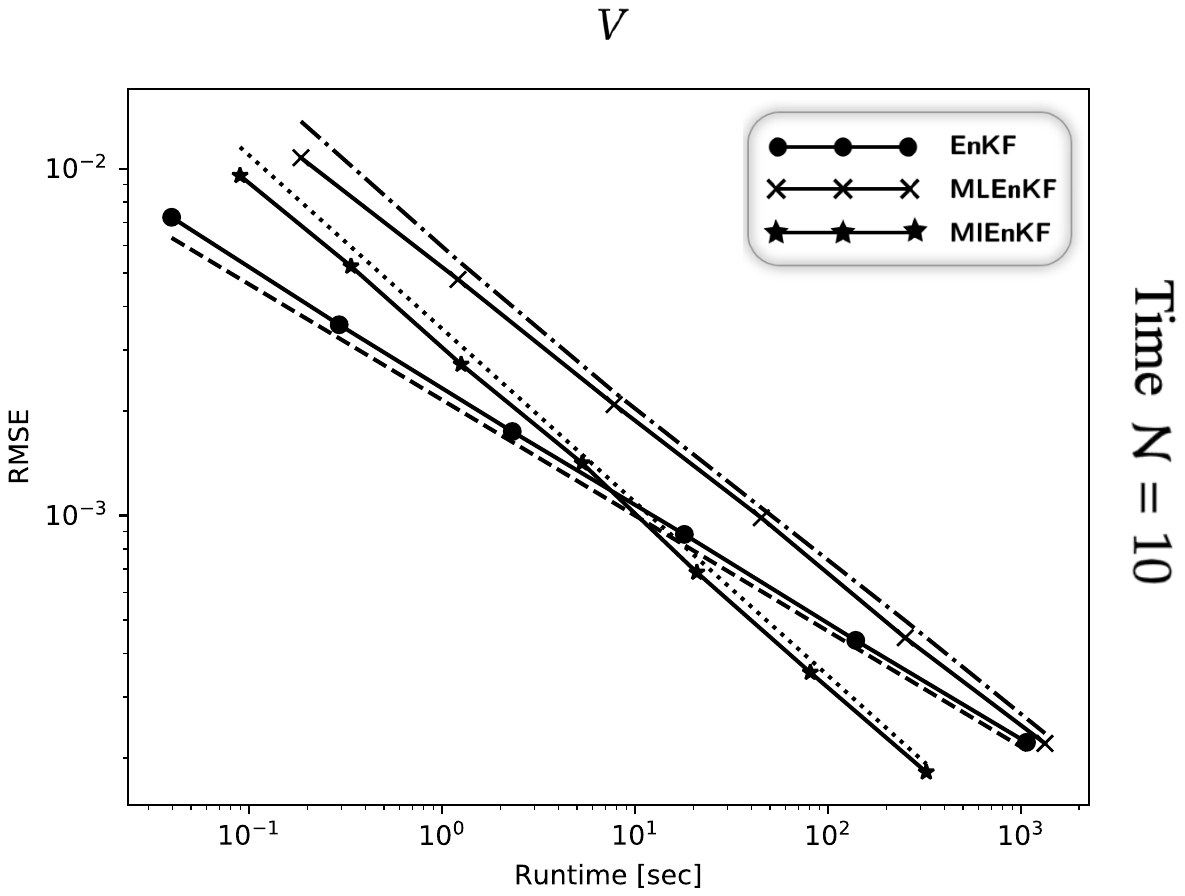}}\\
	\vspace*{-4cm}
		\includegraphics[width=0.47\textwidth]{{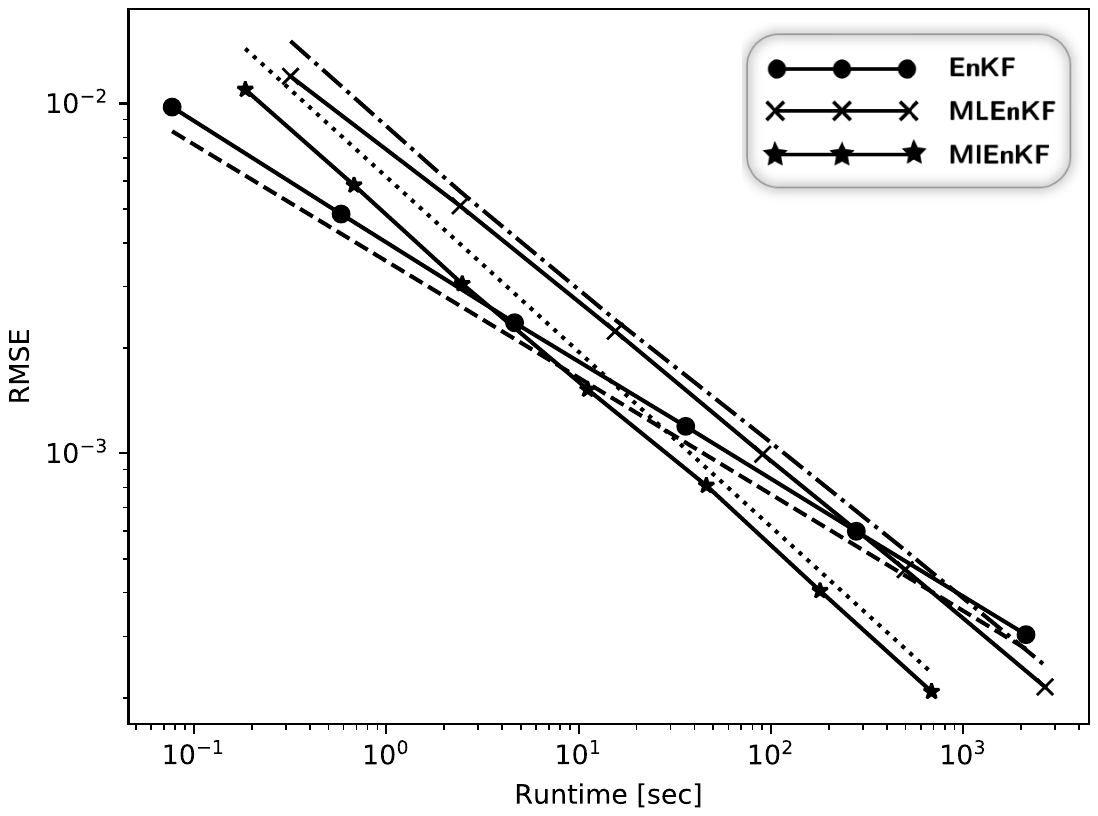}}
	\includegraphics[width=0.47\textwidth]{{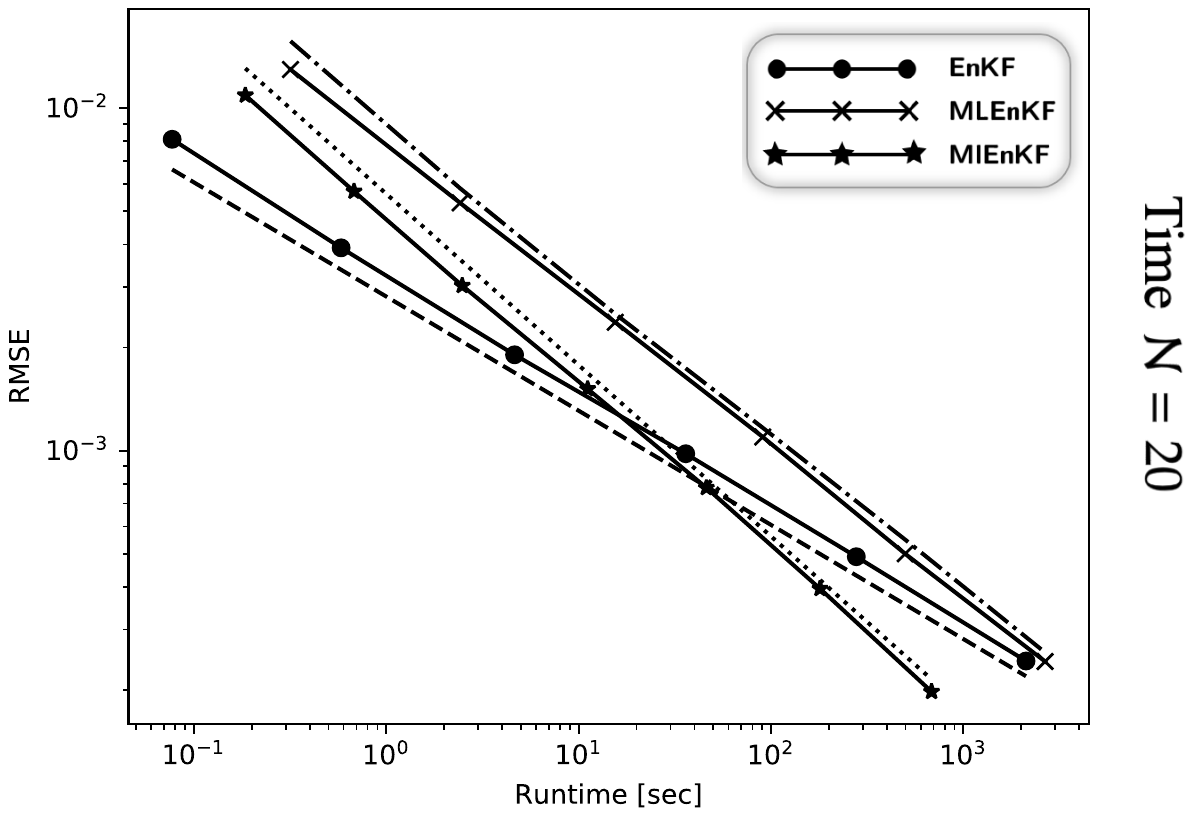}}
		\vspace*{-2cm}
	\caption{\textbf{Langevin dynamics with full observations, $H=[1 \; 0; 0 \; 1]$. Estimates based on $S=90$
			independent runs (Section~\ref{ssec:langevin}).}
	Similar plots as those shown in
	Figure~\ref{fig:loglogRatesLangevinH10sepT10}.}
	\label{fig:loglogRatesLangevinH11sepT10}
\end{figure}

To shed some light on the importance of the temperature parameter,
Figure~\ref{fig:LangevinEvolutionT50Split} illustrates the
phase-portrait time evolution of the realization of Langevin dynamics
up to the final time $\cN=50$ for different temperatures $T=[0, \; 0.01,\;  0.1, \; 1.0]$. Damping
causes a rapid decay of the velocity from the initial value to zero
when $T=0$. For positive temperatures, thermal fluctuation leads to
more diffusive dynamics.
Figure~\ref{fig:LangevinMeanTruthObsT50SplitH11} shows the
signal-tracking performance of MIEnKF for the full observation operator
\[
H=\begin{bmatrix} 1 & 0\\ 0 &1
\end{bmatrix}
\]
and the partial observation operators $H=[1\; 0]$ or $H=[0\; 1]$, all
computed at the tolerance $\epsilon = 2^{-7}$. The method is tracking
the true state of the observed components well in all cases, but, as
is to be expected, it does not track the true state of unobserved
components with the same level of accuracy. The numerical verification
of assumptions~\eqref{ass:A1} and~\eqref{ass:A2star} with respect to
different observation operators is shown in
Figures~\ref{fig:LangevinConjRatesH10}
and~\ref{fig:LangevinConjRatesH01}, respectively. For a sequence of
predefined tolerances, $\epsilon = [2^{-4}, 2^{-5}, \ldots, 2^{-9}]$
for EnKF and MLEnKF and $\epsilon = [2^{-4}, 2^{-5}, \ldots,
2^{-10}]$ for MIEnKF, we compare the performance of the three methods
in terms of runtime versus RMSE. We consider 
$\mathcal{N}=10$ and $\mathcal{N}=20$ observation times, the QoI
$\varphi(X,V)=X$ and $\varphi(X,V)=V$, and we use $S=90$
independent runs of each method to estimate both RMSE and runtime. 
Figures~\ref{fig:loglogRatesLangevinH10sepT10} and~\ref{fig:loglogRatesLangevinH11sepT10}
show the results for the observation operators 
\[
H=[1\; 0] \quad \text{ and } \quad H=\begin{bmatrix} 1 & 0\\ 0 &1
\end{bmatrix},
\]
respectively. The observed complexity
rates for MIEnKF are close to the theory, and the method is more
efficient than the alternatives for small tolerances in both cases.

\section{Conclusion}\label{sec:conclusion}
We have developed a hierarchical ensemble-based filtering method
called the MIEnKF method.  MIEnKF is based on independent samples of
four-coupled EnKF estimators on a multi-index hierarchy of resolution
levels. Under Assumptions~\ref{ass:Psi} and~\ref{ass:Psi2}, we proved
that the method is highly efficient and that it will
asymptotically outperform the comparable methods EnKF and MLEnKF.
For instance, when the weak convergence rate $\alpha=1$ and the
strong convergence rate $\beta=2$, which is a more robust setting of
the EnKF and MLEnKF methods considered
in~\cite{hoel2020multilevel}, the computational cost of reaching
$\cO(\epsilon^2)$ MSE is $\mathcal{O}(\epsilon^{-2})$ for MIEnKF,
$\mathcal{O}(\epsilon^{-2}|\log(\epsilon)|^3)$ for MLEnKF, and
$\mathcal{O}(\epsilon^{-3})$ for EnKF.

In this work we have constructed a multi-index EnKF method with two
resolution parameters: $N_{\ell_1}$ relating to the
time-discretization, and the ensemble-size $P_{\ell_2}$, a 2-index
MIEnKF method. For more complicated high-dimensional
filtering problems, it is an open question if it is possible to extend
MIEnKF to having more resolution parameters, and whether that would
lead to further performance gains. One extension we currently working
on is a 3-index MIEnKF for spatiotemporal models that are discretized
in both space and time, e.g., reaction-diffusion stochastic partial
differential equations (SPDE)~\cite{chernov2020}.

Another interesting direction would be MIEnKF for filtering problems
with high-frequency or continuous-time observations. Here, the new
challenge is that low-resolution levels have to be updated --
has to assimilate observations -- at a lower frequency than high-resolution
levels, but strong coupling still has to be preserved. The recent
work on MLEnKF for Kalman-Bucy filters~\cite{chada2020multilevel}
would be a good starting point for developing an MIEnKF method for
such problem settings.

\medskip

{\bf Acknowledgments } This work was supported by the KAUST Office of
Sponsored Research (OSR) under Award No. URF/1/2584-01-01 and the
Alexander von Humboldt Foundation. G.~Shaimerdenova and R.~Tempone are
members of the KAUST SRI Center for Uncertainty Quantification in
Computational Science and Engineering.

\bibliography{papers}
\bibliographystyle{plain}

\end{document}